\documentclass[final,onefignum,onetabnum]{article}

\usepackage{array}

\usepackage[caption=false]{subfig}
\usepackage{pgfplots}
\usepackage{cite}

\newtheorem{thm}{Theorem}[section]

\newtheorem{rem}[thm]{Remark}

\usepackage{algorithmic}

\usepackage{graphicx,epstopdf}

\usepackage{subfig}
\usepackage{amsmath}
\usepackage{amssymb}
\usepackage{graphicx}
\usepackage{dcolumn}
\usepackage{mathtools}
\usepackage{amsmath,amsfonts}
\usepackage{bm}
\usepackage{amsmath}
\usepackage{amssymb}
\usepackage{color}
\usepackage{float}
\usepackage{setspace}
\usepackage{tabularx}

\allowdisplaybreaks

\newcommand{\f}{\frac}

\newcommand{\E}{ {\mathbb{E}} }

\newcommand{\be}{\begin{equation}}
\newcommand{\ee}{\end{equation}}

\def\lan{\langle} \def\ran{\rangle}

\def\ba{\begin{array}}                \def\ea{\end{array}}
\def\bel{\begin{equation}\label}      \def\ee{\end{equation}}

\colorlet{texcscolor}{blue!50!black}
\colorlet{texemcolor}{red!70!black}
\colorlet{texpreamble}{red!70!black}
\colorlet{codebackground}{black!25!white!25}

\date{}

\title{A Stochastic Maximum Principle Approach for Reinforcement Learning with Parameterized Environment}
\author{
Richard Archibald \thanks{ Division of Computational Science and Mathematics, Oak Ridge National Laboratory.}
\and Feng Bao\thanks{ Department of Mathematics, Florida State University, Tallahassee, Florida, \ ({\tt bao@math.fsu.edu}).} 
\and Jiongmin Yong \thanks{ Department of Mathematics, University of Central Florida, Orlando, Florida.} 
      }

\begin{document}
\maketitle

%% ------------------------------------------------------------------
%% ABSTRACT
%% ------------------------------------------------------------------
%\begin{tcbverbatimwrite}{tmp_\jobname_abstract.tex}
\begin{abstract}
In this work, we introduce a stochastic maximum principle (SMP) approach for solving the reinforcement learning problem with the assumption that the unknowns in the environment can be parameterized based on physics knowledge. For the development of numerical algorithms, we shall apply an effective online parameter estimation method as our exploration technique to estimate the environment parameter during the training procedure, and the exploitation for the optimal policy will be achieved by an efficient backward action learning method for policy improvement under the SMP framework. Numerical experiments will be presented to demonstrate that our SMP approach for reinforcement learning can produce reliable control policy, and the gradient descent type optimization in the SMP solver requires less training episodes compared with the standard dynamic programming principle based methods.

%Different from dynamic programming principle based methods, which solve the original problem via smaller nested sub optimization problems, we can derive a gradient with respect to control actions for policy improvement, which allows us to   Our numerical algorithm tis composed of a backward action learning method for solving the stochastic optimal control problem under the SMP framework and a direct filter method for online estimating the environment parameters.

\end{abstract}


%\textbf{AMS:} 

%\end{tcbverbatimwrite}
%\input{tmp_\jobname_abstract.tex}
%% ------------------------------------------------------------------
%% END HEADER
%% ------------------------------------------------------------------

\section{Introduction}

Reinforcement learning (RL) is an important research area in machine learning. Different from supervised learning and unsupervised learning, RL aims to find how to map situations (state) to actions (control), and the goal of the RL problem is to let an agent learn how take actions in an environment in order to minimize a performance cost or maximize a reward.   As a major machine learning task, RL has been extensively studied and it has application potentials to solve many real-life problems, e.g., robotic automation, natural language processing, health care, image processing, and trading in financial market.  In addition to its straightforward engineering style applications, RL has also drawn increasing attention from the science community. Some recent studies show that RL techniques can be used to solve scientific problems related to dynamic experimental design in physics and chemistry \cite{DRL_ED_2022, ED_RL_21}.

The mathematical foundation of the standard approach for solving the RL problem is the \textit{dynamic programming principle} (DPP) \cite{XYZ_RL_20}, which was introduced to solve the optimal control problem.  The main idea of the DPP approach is to consider a family of local optimal control problems with different initial states and times and establish relationships among these sub-problems through the Hamilton-Jacobi-Bellman equation \cite{Yong_control}. An important numerical method for implementing the DPP for solving the RL problem is temporal difference (TD) learning \cite{TD_Learning_95}, and a major breakthrough was the development of an off-policy TD control algorithm known as Q-learning \cite{Q_learning_92, Q_learning_94}. The Q-Learning method can carry out TD learning efficiently, and it can learn how to take optimal actions without requiring an environment model. This makes Q-Learning applicable to solve many control problems in real-life. On the other hand, TD learning also has some limitations. For example, TD learning methods typically use gradient-free optimization to determine the optimal policy, hence even the Q-learning method suffers from the efficiency issue due to the low convergence rate of gradient-free optimization. Another notable disadvantage of TD learning is that it's a bootstrap procedure that analyzes how good is a guess from another guess. Hence the feedback in TD Learning is delayed and heavily corrupted by noise \cite{Q_Learning_1997fuzzy, Q_learning_94}. As a result, decisions made by TD Learning are more reliable to optimize the agent's short-term performance. However, TD Learning can be misled when short-term gains disagree with the long-term goal. In this case, the agent can be attracted by carefully designed baits that lead it towards a trap. This is especially more challenging when using TD Learning to solve a continuous problem since approximation of the original problem can make predictions more complicated \cite{Q_learning_continuous_16}. Such a drawback in TD Learning is mainly due to the nature of DPP, which solves the optimal control problem by combining a set of local sub-problems, and the fact that no overall physics model is considered to supervise the policy improvement.  

\vspace{0.5em}

In this work, we develop a fundamental methodology of solving the RL problem by using the \textit{stochastic maximum principle} (SMP), which is a major alternative approach for solving the stochastic optimal control problem besides the DPP, and we will focus on the continuous time-space state model with noise perturbations. An important assumption we need for our SMP approach is that the environment is described by a physics model, and the unknown factors in the environment are characterized by model parameters. Although the capability of training an agent how to take actions without an explicit environment model is necessary for many RL applications, being able to take the first principle into account is also important when applying RL to solve scientific problems with well-established physics knowledge. 

As a fundamental mathematical effort, we will not discuss deep learning related deep reinforcement learning techniques (e.g. \cite{DRL_16}), which might be more powerful but certainly will lead to several challenges caused by computational implementation of deep neural networks, such like the overfitting issue, the representation capability, and the reliability in training. Although we don't consider the application of deep learning in this work, appropriately designed deep learning methods can also be applied to our general computational framework and improve the efficiency of our maximum principle approach.
\vspace{0.5em}

For solving the classic stochastic optimal control problem, the SMP aims to find the optimal control by optimizing a stochastic system called the Hamiltonian, and a gradient process with respect to control can be derived by applying the ``G\^ateaux derivative'' to the adjoint process of the state dynamics in order to carry out gradient descent optimization \cite{Peng_control}. The SMP approach has several advantages over the DPP approach. For example, it allows to have random coefficients in the state equation and in the performance cost or reward, and it allows more state constraints -- especially some finite dimensional terminal state constraints (see \cite{Yong_control}). Moreover, the gradient process can give us a direction to improve the policy over the performance period. This is potentially more efficient compared with the DPP approach due to the application of gradient-based optimization. In addition, the gradient with respect to control derived under the SMP framework is based on the understanding of the ``global'' environment.  In this way, the SMP designed policy can better balance between the short-term gains and the long-term goal. Therefore, the SMP approach can be more reliable than the DPP based methods, which typically rely on solutions of stacked local optimization problems. 

\vspace{0.5em}

On the other hand, the SMP approach has to face two challenges when solving the RL problem. First of all, although the gradient process with respect to control can be derived and formulated explicitly by using the adjoint of the state process, obtaining a numerical approximation for the adjoint state process, which is a backward stochastic differential equation, is a challenging task. As a result, the optimization procedure for finding the optimal policy (i.e. the exploitation procedure) needs to evaluate the gradient repeatedly, which makes the SMP approach computationally expensive \cite{GPM_2017}. Secondly, the environment (i.e. the state model and the cost/reward model) in a RL problem is not completely known. Therefore, an environment estimation method for the purpose of exploration is needed while searching for the optimal policy.

%\vspace{1em}

To address the efficiency issue in searching the optimal policy, we introduce a ``backward action learning'' (BAL) method for efficiently solving the stochastic optimal control problem. The main theme of the BAL method is to apply a sample-wise  numerical scheme to approximate the solution of the backward adjoint equation sample-by-sample and then adopt the methodology of stochastic approximation to carry out a stochastic gradient descent procedure to determine the optimal policy. In this way, we can avoid the high computational cost of solving backward stochastic differential equations in the state space. At the same time, the mathematical expression of the gradient process, which contains solutions of the adjoint equation, can still be effectively utilized through the sample-wise approximator to search for the optimal policy.

To address the second issue and explore the environment, we assume that we have enough physics knowledge of the environment so that the environment can be formulated as a parameterized model. Then, exploring the environment is equivalent to searching for the environment parameter, and therefore exploration can be achieved via parameter estimation. In this work, we apply the direct filter, which is an accurate and efficient online parameter estimation method \cite{Bao_parameter}, to estimate the environment parameter during the training procedure. Other online parameter estimation methods may also be used under our general SMP methodology.

Since an online parameter estimation method can dynamically provide feedbacks through training trials and generate real-time updates to improve the understanding of the environment, the estimated environment parameter can guide the BAL optimal control solver to exploit the optimal policy. To further explore the environment and balance between exploration and exploitation, we also perturb the policy by some artificial noise and adopt an ``$\epsilon$-greedy'' mechanism to encourage exploration \cite{epsilon_greedy11}.

\vspace{0.5em}

The rest of this paper is organized as follows. In Section \ref{Methodology}, we introduce the model-based RL problem with parameterized environment, and we shall provide a general methodology of using the SMP to solve the RL problem. In Section \ref{Algorithms}, we introduce the numerical algorithms about how to efficiently solve the RL problem under the SMP framework.  In Section \ref{Numerics}, we present three numerical examples to demonstrate the effectiveness of our algorithm and the necessity of using the SMP approach to solve the RL problem in application scenarios. Some concluding remarks that summarize our research outcomes will be given in Section \ref{Conclusions}.

%\newpage

\section{Problem setting and methodology}\label{Methodology}

In this section, we first introduce the formulation of model-based reinforcement learning (RL) with parameterized environment. Then, we shall introduce a direct filter method for learning the environment parameter and a stochastic maximum principle (SMP) type optimal control solver to find the optimal policy.

\subsection{Problem setting}

In this work, we consider the RL problem under the complete filtered probability space $(\Omega, \mathcal{F}, \mathbb{F}^W, \mathbb{P})$, where $\mathbb{F}^W := \{\mathcal{F}_s^W\}_{s\geq 0}$ is the natural filtration augmented by all the $\mathbb{P}$-null sets in $\mathcal{F}$. The dynamics of the agent is formulated by the following stochastic dynamical system in the form of a stochastic differential equation (SDE)
\begin{equation}\label{Control:State}
dX^{\bm{a}, \lambda}_t = b(t, X^{\bm{a}, \lambda}_t, \bm{a}_t, \lambda) dt + \sigma(t, X^{\lambda}_t, \bm{a}_t, \lambda) dW_t, \hspace{1em} t \in [0, T],
\end{equation}
where $X^{\bm{a}, \lambda}_t \in \mathbb{R}^d$ is the state of the agent, which is also known as the \textit{``state process''} (or the ``state equation'') in the optimal control problem; $W:= \{W_t\}_{t\in [0, T]}$ is a standard $d$-dimensional Brownian motion that introduces uncertainty to the agent state; $b: [0, T] \times \mathbb{R}^d \times \mathbb{R}^m \times \mathbb{R}^q \rightarrow \mathbb{R}^d$ and $\sigma: [0, T] \times \mathbb{R}^d \times \mathbb{R}^m \times \mathbb{R}^q \rightarrow \mathbb{R}^{d \times d}$ are suitable maps called \textit{drift} and \textit{diffusion}, respectively. In the model-based RL problem, we use $b$ and $\sigma$ to model the dynamics of the agent, and they contain the ``physics knowledge'' of the environment that we already possess. Since we assume that there are also unknowns in the environment, we introduce a state-dependent parameter $\lambda: \mathbb{R}^d \rightarrow \mathbb{R}^q$ to represent \textit{physics-informed} unknown factors that we need to learn in the environment. The vector-valued process $\bm{a}_t \in \mathbb{R}^m$ stands for the agent action at time $t$, which is equivalent to the control process in the stochastic optimal control problem. In the RL problem, the control process is also called the ``policy''. Denote $\mathcal{U}[0, T] : = \Big\{\bm{a}: [0, T] \times \Omega \rightarrow U \subset \mathbb{R}^m \big| \text{ $\bm{a}$ is $\mathbb{F}^W$-progressively measurable} \Big\}$ as the admissible control set, in which we can choose the control actions.
Under some mild conditions \cite{Yong_control}, for every choice of parameter $\lambda \in \mathbb{R}^q$ and control $\bm{a} \in \mathcal{U}[0, T]$, SDE \eqref{Control:State} admits a unique solution $X^{\bm{a}, \lambda}$. 

The performance of the action $\bm{a}$ is measured by the following cost (or reward):
\begin{equation}\label{cost}
J(\bm{a}) = \E\left[ \int_0^T f^{\lambda}(t, X^{\bm{a}, \lambda}_t, \bm{a}_t) dt + h(X^{\bm{a}, \lambda}_T) \right],
\end{equation}
where $f^{\lambda}$ is the running cost and $h$ measures the cost at the terminal time $T$. 

In this work, we let the ``environment'' in the RL problem be mathematically formulated by $b$, $\sigma$, $f^{\lambda}$ and $h$. The unknowns in the environment are represented by parameter $\lambda$, and we want to re-emphasize that being able to incorporate physics knowledge into a RL task is necessary in many practical scientific machine learning problems.

The goal of the stochastic optimal control problem is to find the ``\textit{optimal control}'' $\bar{\bm{a}}$ that minimizes the cost $J$ \footnote{In the case of $J$ is a reward, we maximize the reward functional.}, i.e.
\begin{equation}\label{minimize}
J(\bar{\bm{a}}) = \inf_{\bm{a} \in \mathcal{U}[0, T]} J(\bm{a}).
\end{equation}
In the RL language, we aim to find the optimal policy $\bar{\bm{a}}$ that minimizes the ``penalty''. Similar framework can also search for the optimal policy that maximizes the ``reward'' by switching the minimization problem to maximization problem in Eq. \eqref{minimize}. 

When parameter $\lambda$ is given, equations \eqref{Control:State} - \eqref{minimize} formulate a classic stochastic optimal control problem. The major difference that distinguishes the RL problem from the stochastic optimal control problem is the unknown environment, which is modeled by the physics-informed unknown parameter $\lambda$ in this work. 
\vspace{0.75em}

As a numerical approach for solving the RL problem, our method will also consist an exploration procedure and an exploitation procedure. Since the environment in this work is parameterized, the exploration procedure, which aims to determine the environment during learning, is equivalent to implement online parameter estimation for $\lambda$. On the other hand, the exploitation procedure, which finds the optimal action, is equivalent to solving the stochastic optimal control problem. In what follows, we shall introduce a SMP type RL solver with online parameter estimation for learning the environment.

\subsection{A stochastic maximum principle framework for the reinforcement learning problem}

Since the major difference between the RL problem and the classic optimal control problem is the exploration of the environment, finding the parameter that represents the environment while learning the optimal policy is a key challenge. In this work, we formulate the dynamical estimation of the environment parameter as an optimal filtering problem.

\subsubsection*{Exploration: Learning the environment by using the direct filter method}

The main idea of the ``optimal filtering'' approach for dynamically estimating parameters is to ``project'' the data of the agent state to the parameter space and use the conditional probability density function (PDF) of the parameter, i.e. $p(\lambda_k | \mathcal{X}_k)$, to calculate the estimate for $\lambda$ at the $k$-th episode, where $\sigma$-algebra $\mathcal{X}_k : = \sigma\big(X^{\bm{a}_l, \lambda_l} \vee f^{\lambda_l}, 0 \leq t \leq T, l = 0, 1, \cdots, k\big)$ contains the information of the state of the agent as well as the parameterized cost in the previous training episodes. In this work, we apply the direct filter method \cite{Bao_parameter} to estimate the environment parameter in the online manner.

To proceed, we use a sequence of random variables $\{\lambda_k\}_k$ to represent our estimates for the unknown parameter $\lambda$, where $\lambda_k$ is the estimate corresponding to the $k$-th training episode. Assuming that we have $\lambda_k$ that follows the conditional PDF $p(\lambda_k | \mathcal{X}_k)$, we generate a proposal parameter random variable, i.e. $\lambda_{k+1}$, through the following pseudo-dynamics:
\begin{equation}\label{dynamics:lambda}
\lambda_{k+1}= \lambda_{k} + \xi_k, \hspace{2em} k = 0, 1, 2, \cdots,
\end{equation}
where $\xi_k \sim N(0, (\delta_{k})^2)$ is a standard Brownian motion with pre-chosen covariance constants $\{\delta_{k}\}_{k\geq 0}$, and $\xi_k$ gives artificial noise that allows to explore other possible values of the environment parameter.

For a given conditional PDF $p(\lambda_k | \mathcal{X}_k)$, which describes the estimated parameter at the $k$-th episode, we apply the following Bayesian inference to obtain the posterior PDF $p(\lambda_{k+1} | \mathcal{X}_{k+1})$ that ``optimally'' describes the estimated parameter $\lambda_{k+1}$ corresponding to the state information of the $k+1$-th episode as follows
\begin{equation}\label{Bayesian:lambda}
p(\lambda_{k+1} | \mathcal{X}_{k+1}) \sim p(\lambda_{k+1} | \mathcal{X}_k) p(X^{\bm{a}, \lambda}, f^{\lambda} | \lambda_{k+1}),
\end{equation}
where $p(\lambda_{k+1} | \mathcal{X}_k)$ is the prior conditional PDF of $\lambda_{k+1}$ derived from the pseudo-dynamics \eqref{dynamics:lambda} and the parameter distribution $p(\lambda_k | \mathcal{X}_k)$ at the training episode $k$, and $p(X^{\bm{a}, \lambda}, f^{\lambda} | \lambda_{k+1})$ is the likelihood function that compares the simulated agent state and its corresponding cost derived from the prior parameter variable $\lambda_{k+1}$ with the real agent state $X^{\bm{a}, \lambda}$ and the real cost $f^{\lambda}$ produced by the true environment parameter $\lambda$.

\subsubsection*{Exploitation: The stochastic maximum principle approach for searching the optimal policy}

In the case that the complete knowledge of the environment ( or the environment parameter $\lambda$ is known), and the optimal strategy $\bar{\bm{a}}$ is in the interior of $U$ \footnote{In many practical RL scenarios, people often let $U$ be the real space.}, we can deduce by using the G\^ateaux derivative of $\bar{\bm{a}}$ and the maximum principle that the gradient process of the cost functional $J$ with respect to the control process over time interval $t \in [0, T]$ has the following form \cite{Yong_control, GPM_2017}
\begin{equation}\label{Gradient-Standard}
\nabla J_{\bm{a}} (\bar{\bm{a}}_t) = \E\left[b_{\bm{a}}(t, \bar{X}^{\bar{\bm{a}}, \lambda}_t, \bar{\bm{a}}_t, \lambda)^{\top} \bar{Y}^{\bar{\bm{a}}, \lambda}_t + \sigma_{\bm{a}}(t, \bar{X}^{\bar{\bm{a}}, \lambda}_t, \bar{\bm{a}}_t, \lambda)^{\top} \bar{Z}^{\bar{\bm{a}}, \lambda}_t  + f^{\lambda}_{\bm{a}}(t,  \bar{X}^{\bar{\bm{a}}, \lambda}_t, \bar{\bm{a}}_t)^{\top} \right],
\end{equation}
where subscripts are used to denote partial derivatives of functions. The stochastic processes $\bar{Y}$ and $\bar{Z}$ are adapted solutions of the following forward backward stochastic differential equations (FBSDEs) system
\begin{equation}\label{FBSDEs}
\begin{aligned}
d\bar{X}^{\bar{\bm{a}}, \lambda}_t &= b(t, \bar{X}^{\bar{\bm{a}}, \lambda}_t, \bar{\bm{a}}_t, \lambda) dt + \sigma(t, \bar{X}^{\bar{\bm{a}}, \lambda}_t, \bar{\bm{a}}_t, \lambda) dW_t, \qquad & \text{(forward SDE)}\\
d\bar{Y}^{\bar{\bm{a}}, \lambda}_t &=  \big( - b_x(t, \bar{X}^{\bar{\bm{a}}, \lambda}_t, \bar{\bm{a}}_t, \lambda)^{\top} \bar{Y}^{\bar{\bm{a}}, \lambda}_t - \sigma_x(t, \bar{X}^{\bar{\bm{a}}, \lambda}_t, \bar{\bm{a}}_t, \lambda)^{\top} \bar{Z}^{\bar{\bm{a}}, \lambda}_t - f^{\lambda}_x(t, \bar{X}^{\bar{\bm{a}}, \lambda}_t, \bar{\bm{a}}_t)^{\top} \big) dt \\
 & \qquad + \bar{Z}^{\bar{\bm{a}}, \lambda}_t dW_t, \qquad \bar{Y}^{\bar{\bm{a}}, \lambda}_T = h_x(\bar{X}^{\bar{\bm{a}}, \lambda}_T)^{\top}, & \text{(BSDE)}
\end{aligned}
\end{equation}
where the first equation in \eqref{FBSDEs} is a standard forward stochastic differential equation (SDE) with the same expression as the state equation for the agent, and the second equation is a backward stochastic differential equation (BSDE), which is also call the ``adjoint equation'' of the state equation \eqref{Control:State}. The solution $Z$ of the BSDE is the martingale representation of $Y$ with respect to the Brownian motion $W$.

With gradient $\nabla J_{\bm{a}}$ introduced in Eq. \eqref{Gradient-Standard} and solutions $X$, $Y$, and $Z$ introduced in the FBSDE system \eqref{FBSDEs}, the SMP approach for solving the stochastic optimal control problem will carry out the following gradient descent optimization procedure to solve for the optimal control, i.e. the optimal policy in RL,
\begin{equation}\label{Gradient-Descent}
\bm{a}^{k+1} = \bm{a}^{k} - \eta_k \nabla J_{\bm{a}}(\bm{a}^{k}), \qquad k = 0, 1, 2, \cdots, K-1,
\end{equation}
where we have an initial guess policy $\bm{a}^{0}$, $\eta_k$ is the learning rate, $K$ is a pre-determined total number of training episodes, and we let $\bm{a}^{K}$ be our estimated optimal policy.

\subsubsection*{The general methodology of solving the reinforcement learning problem via the stochastic maximum principle}
In what follows, we briefly discuss the general methodology about how to solve the RL problem under the SMP framework. At this moment, we ignore computational implementation issues, and the numerical algorithms will be introduced in the next section.

The main theme of the SMP solver for RL is to use the iterative scheme \eqref{Gradient-Descent} to improve the policy and use the direct filter to explore the environment.  Assume that with $k$ episodes of training, we have an estimated environment parameter $\lambda_k$ (with its conditional PDE $p(\lambda_{k}|\mathcal{X}_k)$) and a sub-optimal policy $\bm{a}^{k}$ based on the understanding of the environment corresponding to $\lambda_k$. 

For the $k+1$-th episode, we use scheme \eqref{dynamics:lambda} to generate a prior conditional PDF $p(\lambda_{k+1}|\mathcal{X}_k)$ that characterizes the proposal parameter. Then, we let the agent follow the current policy $\bm{a}^{k}$ and interact with the real environment to obtain a state process $X^{\bm{a}^{k}, \lambda}$, which can be used to compare with simulated state processes to generate the likelihood $p(X^{\bm{a}^k, \lambda}_{t} | \lambda_{k+1})$, and we can obtain the posterior PDF $p(\lambda_{k+1} | \mathcal{X}_{k+1})$ for the estimated parameter variable $\lambda_{k+1}$ via the Bayesian inference \eqref{Bayesian:lambda}. Note that the agent trial state process, i.e. $X^{\bm{a}^{k}, \lambda}$, provides the feedback from the real environment for exploration.  With an updated estimate for the environment at the learning episode $k+1$, we apply the SMP method to find the optimal policy. Specifically, we compute the solutions of the FBSDE system \eqref{FBSDEs} with the estimated policy $\bm{a}^{k}$ and the environment parameter $\lambda_{k+1}$ to obtain approximations  for $X^{\bm{a}^{k}, \lambda_{k+1}}$, $Y^{\bm{a}^{k}, \lambda_{k+1}}$ and $Z^{\bm{a}^{k}, \lambda_{k+1}}$.   Then, the gradient $\nabla J_{\bm{a}}$ introduced in Eq. \eqref{Gradient-Standard} can be calculated by using the simulated solutions of the FBSDE system \eqref{FBSDEs}. As a result, the approximated gradient would give us a direction to improve the current policy, and we can obtain the improved policy $\bm{a}^{k+1}$ via the gradient descent scheme \eqref{Gradient-Descent}. To encourage exploration, we also adopt the ``$\epsilon$-greedy'' method by perturbing the sub-optimal policy with some artificial noise.

\section{Numerical algorithms}\label{Algorithms}

In this section, we derive numerical algorithms to implement the above SMP approach for the RL problem.  For convenience of presentation, we assume that the diffusion coefficient $\sigma$ in the state dynamics \eqref{Control:State} is a deterministic time-dependent process, denoted by $\sigma_t$.  Algorithms for more general cases can be obtained under our methodology with more tedious derivation, and the model with the simplified diffusion coefficient can already cover wide range of application problems since the physics knowledge of a stochastic model is often incorporated into the drift term. In what follows, we shall first discuss numerical implementation of the direct filter method for exploration (Section \ref{Direct-Filter}), and we will provide a \textit{backward action learning (BAL)} method for solving the classic stochastic optimal control problem under the SMP framework for the purpose of exploitation (Section \ref{BAL-SOC}). Then, we combine exploration with exploitation and introduce an overarching algorithm to solve the RL problem (Section \ref{BAL-RL}).   

\subsection{Particle implementation of the direct filter method for exploration}\label{Direct-Filter}
In this paper, we adopt the numerical recipe of the ``particle filter'', which is also known as a ``sequential Monte Carlo method'', to implement the direct filter for exploring the environment in the RL problem \cite{Bao_Cogan20, particle-filter}. 

\vspace{0.5em}

Assume that after $k$ training episodes, we have a (suboptimal) policy $\bm{a}^k$ and a set of $Q$ samples (called ``particles'' in the particle filter), denoted by $\{\zeta_k^{(q)}\}_{q=1}^Q$, that follow the conditional PDF $p(\lambda_k | \mathcal{X}_k)$ for the estimated parameter $\lambda_k$. To carry out the exploration procedure in the $k+1$-th episode, we first generate a set of particles based on the pseudo-dynamics introduced in Eq. \eqref{dynamics:lambda} to generate a set of proposal particles $\{\tilde{\zeta}_{k+1}^{(q)}\}_{q=1}^Q$ as follows
\begin{equation}\label{parameter:prediction}
\tilde{\zeta}_{k+1}^{(q)} =\zeta_{k}^{(q)} + \xi_k^{(q)}, \qquad q = 1, 2, \cdots, Q,
\end{equation}
where $\xi_k^{(q)} \sim \xi_k$, and the particle set $\{\tilde{\zeta}_{k+1}^{(q)}\}_{q=1}^Q$ formulates the empirical distribution, i.e.
$\tilde{p}(\lambda_{k+1} | \mathcal{X}_k) : = \f{1}{Q} \sum_{q=1}^Q \delta_{\tilde{\zeta}^{(q)}_{k+1}} (\lambda_{k+1})$, 
for the prior  PDF $p(\lambda_{k+1} | \mathcal{X}_k)$.

\vspace{0.5em}

To incorporate the trial agent state at the $k+1$-th training episode and learn the environment, we update the proposal particles through Bayesian inference. Specifically, for each proposal particle $\tilde{\zeta}_{k+1}^{(q)}$, we generate a simulated state trajectory $X^{\bm{a}^k, \tilde{\zeta}_{k+1}^{(q)}} : = \{X_t^{\bm{a}^k, \tilde{\zeta}_{k+1}^{(q)}}\}_{0 \leq t \leq T}$ based on the current policy $\bm{a}^k$ and the proposal particle $\tilde{\zeta}_{k+1}^{(q)}$. On the other hand, the agent that follows policy $\bm{a}^{k}$ interacts with the real environment, and it generates the real state trajectory  $X^{\bm{a}^k, \lambda}$ that reflects the true environment parameter $\lambda$. Then, by comparing each simulated state sample trajectory $X^{\bm{a}^k, \tilde{\zeta}_{k+1}^{(q)}}$ with the real agent state $X^{\bm{a}^k, \lambda}$, we have the following (unnormalized) likelihood for the parameter particle $\tilde{\zeta}_{k+1}^{q}$ as
\begin{equation*}
\begin{aligned}
p(X^{\bm{a}^k, \lambda}, f^{\lambda} \big| \tilde{\zeta}_{k+1}^{(q)}) \sim  \exp\Big( - \Big[ & \big(\int_0^T X_t^{\bm{a}^k, \tilde{\zeta}_{k+1}^{(q)}} dt -  \int_0^T X_t^{\bm{a}^k,\lambda}  dt\big)^2 \\
& \quad + \big(\int_0^T f^{\tilde{\zeta}_{k+1}^{(q)}} (t, X_t^{\bm{a}^k, \tilde{\zeta}_{k+1}^{(q)}} , \bm{a}^k_t)dt -  \int_0^T  f^{\lambda}(t, X_t^{\bm{a}^k, \lambda} , \bm{a}^k_t) dt\big)^2 \Big]/ 2 \delta_k^2\Big),
\end{aligned}
\end{equation*}
where $\delta_k$ (introduced in Eq. \eqref{dynamics:lambda}) is the standard deviation of the artificial noise in the pseudo parameter process that encourages exploration.

Then, combining the prior with the above likelihood through Bayesian inference, we have 
\begin{equation}\label{Bayesian-continuous}
p\big(\lambda_{k+1} = \tilde{\zeta}_{k+1}^{(q)} \big| \mathcal{X}_{k+1} \big) = \f{\tilde{p}(\lambda_{k+1} \big| \mathcal{X}_k)  \ p(X^{\bm{a}^k, \lambda}, f^{\lambda} \big| \tilde{\zeta}_{k+1}^{(q)})}{C}, \qquad q = 1, 2, \cdots, Q,
\end{equation} 
where $C$ is a normalization factor. 

Since the empirical prior distribution $\tilde{p}(\lambda_{k+1} \big| \mathcal{X}_k)$ is described by a set of unweighted prediction particles, the weighted particle pairs $\{(\tilde{\zeta}_{k+1}^{(q)}, \omega_{k+1}^{(q)})\}_{q=1}^Q$ obtained in Eq. \eqref{Bayesian-continuous} can describe the (weighted) empirical distribution for $p\big(\lambda_{k+1} \big| \mathcal{X}_{k+1} \big)$, where the weight $\omega_{k+1}^{(q)} := p(X^{\bm{a}^k, \lambda}, f^{\lambda} | \tilde{\zeta}_{k+1}^{(q)})/C$ is the likelihood of each particle. To improve the stability and address the degeneracy of the particles \cite{MTAC2012, cd2002, Kang-PF}, i.e. only a very small number of particles have significant likelihood weights \cite{MCMC-PF}, we also introduce a classic bootstrap resampling procedure by using the importance sampling method \cite{particle-filter} for the weighted particle pairs $\{\tilde{\zeta}_{k+1}^{(q)}, \omega_{k+1}^{(q)}\}_{q=1}^Q$ and obtain a set of equally weighted particles $\{\zeta_{k+1}^{(q)}\}_{q=1}^Q$, which follow the conditional PDF $p\big(\lambda_{k+1} \big| \mathcal{X}_{k+1} \big)$ as needed for the next exploration procedure in the next training episode.

\begin{rem}
Note that the simulated state trajectories $\{ X^{\bm{a}^k, \tilde{\zeta}_{k+1}^{(q)}} \}_{q=1}^Q$ need to be calculated on discrete temporal points with appropriate numerical schemes. Since the state of the agent coincides with the forward SDE in the FBSDE system, we shall postpone our discussion on the numerical method for the state dynamics in the next subsection when we introduce the numerical method for solving FBSDEs.
\end{rem}

\begin{rem}
In the case that the environment parameter $\lambda$ is state-dependent, we let $\lambda$ be a vector corresponding to agent states and carry out the direct filter method to estimate the parameter if the agent enters the correspondent state block. 
\end{rem}

\subsection{Backward action learning for exploitation}\label{BAL-SOC}

To introduce the numerical algorithm for exploitation, which is equivalent to solving a stochastic optimal control problem, we first assume that we have complete knowledge of the environment, i.e. the environment parameter $\lambda$ is known. Then, we shall combine the direct filter method for exploration with the backward action learning method for exploitation to solve the RL problem. 

The computational framework of our backward action learning method is to carry out a gradient descent optimization procedure with iterative scheme \eqref{Gradient-Descent} to improve the policy, and the gradient with respect to the policy, which is introduced in Eq. \eqref{Gradient-Standard}, is derived based on the SMP with usage of the G\^ateaux derivative. Since the gradient is composed of solutions of the FBSDE system, numerical methods for solving FBSDEs are needed.

%The primary computational cost of such a stochastic maximum principle type approach includes simulation of the expectation for the gradient process and approximation for the numerical solution of the FBSDE system \eqref{FBSDEs}. 

Numerical solvers for both forward SDEs and backward SDEs have been well studied \cite{SDE1, Bao_first, Zhao_multi}. In what follows, we introduce the standard numerical schemes for solving SDEs and BSDEs.

\vspace{0.5em}

To proceed, we introduce a temporal partition 
$$\Pi_{N_T} := \{t_n: 0 = t_0 < t_1 < \cdots < t_n < \cdots < t_{N_T} = T\},$$ 
and we consider the following FBSDE system \eqref{FBSDEs} over the time interval $[t_n, t_{n+1}]$ with a (suboptimal) policy $\bm{a}^{k}$ and a given environment parameter $\lambda$, which will be estimated via the direct filter method in the RL problem, 
\begin{equation}\label{FBSDEs:interval}
\begin{aligned}
X^{\bm{a}^k, \lambda}_{t_{n+1}} &= X^{\bm{a}^k, \lambda}_{t_n} + \int_{t_n}^{t_{n+1}}b(t, X^{\bm{a}^k, \lambda}_{t}, \bm{a}_t^k, \lambda) dt + \int_{t_n}^{t_{n+1}}\sigma_t dW_t,  & \text{(forward SDE)}\\
Y^{\bm{a}^k, \lambda}_{t_n} &= Y^{\bm{a}^k, \lambda}_{t_{n+1}} + \int_{t_n}^{t_{n+1}}\Big( b_x(t, X^{\bm{a}^k, \lambda}_t, \bm{a}^{k}_t, \lambda)^{\top} Y^{\bm{a}^k, \lambda}_t + f^{\lambda}_x(t, X^{\bm{a}^k, \lambda}_t, \bm{a}^{k}_t)^{\top} \Big) dt \\
 & \qquad - \int_{t_n}^{t_{n+1}}Z^{\bm{a}^k, \lambda}_{t}  dW_t.  & \text{(BSDE)}
 \end{aligned}
\end{equation}

For a random variable $X_{t_n}^{\bm{a}^k, \lambda}$ that represents the solution of the forward SDE at time $t_n$, we can approximate $X_{t_{n+1}}^{\bm{a}^k, \lambda}$ by using the following standard Euler-Maruyama scheme
\begin{equation}\label{Approx:X}
X_{t_{n+1}}^{\bm{a}^k, \lambda} \approx X_{t_{n}}^{\bm{a}^k, \lambda} + b(t_n, X^{\bm{a}^k, \lambda}_{t_n}, \bm{a}_{t_n}^k, \lambda) \Delta t_n + \sigma_{t_n} \Delta W_{t_n},
\end{equation}
where $\Delta t_n := t_{n+1} - t_n$ and $\Delta W_{t_n} := W_{t_{n+1}} - W_{t_n}$.

\vspace{0.5em}

To solve the BSDE and obtain a numerical solution for $Y$, we take conditional expectation $\E_n^{X^k}[\cdot] : = \E[\cdot | X_{t_n}^{\bm{a}^k, \lambda}]$ on both sides of the BSDE in Eq. \eqref{FBSDEs:interval}. Since the BSDE is the adjoint equation of the forward state equation, which is backward in time, we assume that a random variable $Y_{t_{n+1}}^{\bm{a}^k, \lambda}$ representing the solution of the BSDE at time $t_{n+1}$ is given, and we use the right-point formula to approximate the deterministic integral on the right hand side of the BSDE. Then, we obtain the following approximation scheme for $Y$ \cite{Bao_IJUQ19}
\begin{equation}\label{Approx:Y}
Y_{t_n}^{\bm{a}^k, \lambda} \approx  \E_n^{X^k}\big[ Y_{t_{n+1}}^{\bm{a}^k, \lambda}\big] + \E_n^{X^k}\Big[ b_x(t_{n+1}, X^{\bm{a}^k, \lambda}_{t_{n+1}}, \bm{a}^{k}_{t_{n+1}}, \lambda)^{\top} Y^{\bm{a}^k, \lambda}_{t_{n+1}} + f^{\lambda}_x(t_{n+1}, X^{\bm{a}^k, \lambda}_{t_{n+1}}, \bm{a}^{k}_{t_{n+1}})^{\top} \Big] \Delta t_n,
\end{equation}
where we have used the martingale property of It\^o type stochastic integrals to get $\E_n^{X^k}\left[ \int_{t_n}^{t_{n+1}}Z^{\bm{a}^k, \lambda}_{t}  dW_t\right] = 0$, and note that $Y_{t_n}^{\bm{a}^k, \lambda} = \E_n^{X^{k}}[Y_{t_n}^{\bm{a}^k, \lambda}]$ due to the adaptedness of $Y$ with respect to $X$.

\vspace{0.5em}

By using approximation equations \eqref{Approx:X}-\eqref{Approx:Y} as a guideline, we introduce the following (temporal-discretized) scheme for solving the FBSDE system \cite{ZhangJ_BSDE, Bao_Control_20, BSDE_filter}:
\begin{equation}\label{Scheme-Full:FBSDEs}
\begin{aligned}
X_{n+1}^{\bm{a}^k, \lambda} =& X_{n}^{\bm{a}^k, \lambda} + b(t_n, X^{\bm{a}^k, \lambda}_{n}, \bm{a}_{t_n}^k, \lambda) \Delta t_n + \sigma_{t_n} \Delta W_{t_n},\\
Y_{n}^{\bm{a}^k, \lambda} = &  \E_n^{X^k}\big[ Y_{n+1}^{\bm{a}^k, \lambda}\big] +\E_n^{X^k} \Big[ b_x(t_{n+1}, Y^{\bm{a}^k, \lambda}_{n+1}, \bm{a}^{k}_{t_{n+1}}, \lambda)^{\top} X^{\bm{a}^k, \lambda}_{n+1} + f^{\lambda}_x(t_{n+1}, X^{\bm{a}^k, \lambda}_{n+1}, \bm{a}^{k}_{t_{n+1}})^{\top} \Big] \Delta t_n,
\end{aligned}
\end{equation}
where $X_{n+1}^{\bm{a}^k, \lambda}$ and $Y_{n}^{\bm{a}^k, \lambda}$ are approximations for $X_{t_{n+1}}^{\bm{a}^k, \lambda}$ and $Y_{t_n}^{\bm{a}^k, \lambda}$ with an estimated policy $\bm{a}^{k}$ and an environment parameter $\lambda$. The side condition of the BSDE is $Y_{N_T}^{\bm{a}^k, \lambda} = h_x$, where $h$ is the terminal cost (penalty), and $X$ is initialized with the initial state of the agent. With scheme \eqref{Scheme-Full:FBSDEs}, we use numerical solutions $X_{n}^{\bm{a}^k, \lambda}$ and $Y_{n}^{\bm{a}^k, \lambda}$ to approximate $X$ and $Y$ in the gradient process and get the following approximation scheme for $\nabla J_{\bm{a}}$:
\begin{equation}\label{Approx:Gradient}
\nabla J_{\bm{a}} (\bm{a}^{k}_{t_n}) \approx \nabla \hat{J}_{\bm{a}} (\bm{a}^{k}_{t_n})= \E\left[b_{\bm{a}}(t_n, X_{n}^{\bm{a}^k, \lambda}, \bm{a}^{k}_{t_n}, \lambda)^{\top} Y_{n}^{\bm{a}^k, \lambda} + f^{\lambda}_{\bm{a}}(t_n, X_{n}^{\bm{a}^k, \lambda}, \bm{a}^{k}_{t_n})^{\top} \right], \ n = 0,1, \cdots, N_T - 1.
\end{equation}
Then, the iterative scheme for finding the optimal control, i.e. the optimal policy in the RL problem, becomes
\begin{equation}\label{Approx:GD}
\bm{a}_{t_n}^{k+1} = \bm{a}_{t_n}^{k} - \eta_k \nabla \hat{J}_{\bm{a}} (\bm{a}_{t_n}^{k}), \qquad n = 0, 1, 2, \cdots, N_T -1.
\end{equation}

\vspace{0.5em}
To implement the gradient descent iteration \eqref{Approx:GD} with the gradient fully calculated following the approximation scheme \eqref{Approx:Gradient}, one needs to evaluate the expectation $\E[\cdot]$ in Eq. \eqref{Approx:Gradient} as well as the conditional expectation $\E_n^{X^k}[\cdot]$ in scheme \eqref{Scheme-Full:FBSDEs}, which is needed to approximate the numerical solutions $X_{n+1}^{\bm{a}^k, \lambda}$ and $Y_{n}^{\bm{a}^k, \lambda}$. The standard approach to evaluate an expectation (especially in high-dimensional spaces) is the Monte Carlo method, in which we use simulated Monte Carlo samples to represent the random variables and use the average of Monte Carlo samples as an approximation for the desired expected value. However, when utilizing the Monte Carlo method in gradient descent optimization with numerical solution $Y_{n}^{\bm{a}^k, \lambda}$ of the adjoint BSDE, in addition to simulating the expectation for the expected gradient $\nabla J_{\bm{a}}$ with Monte Carlo samples, one also needs to generate a large number of state samples for $X$ at each time step $t_n$ in order to evaluate the conditional expectation $\E_n^{X^k}$ in the numerical scheme \eqref{Scheme-Full:FBSDEs} for the calculation of $Y_{n}^{\bm{a}^k, \lambda}$. Since $Y_t^{\bm{a}, \lambda}$ is a random variable whose value is corresponding to the state $X_t^{\bm{a}, \lambda}$, which is also a random variable that continuously takes values in the state space $\mathbb{R}^d$, a Monte Carlo type representation of $Y$ with a set of random samples requires numerical approximation for $Y$ in the entire state space. This would cause a very challenging computational task of high-dimensional approximation when the state dimension $d$ is large, which is often computationally prohibitive due to the so-called ``curse of dimensionality''.  Moreover, the numerical approximation for expectation $\E$ and conditional expectation $\E_n^{X^k}$ needs to be calculated repeatedly over the gradient descent iteration procedure, which will make the full calculation of the gradient descent optimization procedure infeasible in practice.
 
\vspace{0.5em}

In our backward action learning approach for the stochastic optimal control problem, which is inspired by the application of stochastic approximation in gradient descent optimization, we use a single realization of Monte Carlo sample (or a mini-batch of samples) in the Monte Carlo approximation to represent the entire state of the random variable in an expectation. Specifically, at each iteration stage $k$ we use the Euler-Maruyama scheme to generate \textit{one} sample-path $\{\tilde{X}_{n+1}^{\bm{a}^k, \lambda}\}_{n=1}^{N_T}$ for the state process as follows
\begin{equation}\label{SDE-sample}
\tilde{X}_{n+1}^{\bm{a}^k, \lambda} = \tilde{X}_{n}^{\bm{a}^k, \lambda} + b(t_n, \tilde{X}^{\bm{a}^k, \lambda}_{n}, \bm{a}_{t_n}^k, \lambda) \Delta t_n + \sigma_{t_n} \sqrt{\Delta t_n} \epsilon_n, \quad n = 0, 1, 2, \cdots, N_T - 1,
\end{equation}
where $\epsilon_n$ is a random sample drawn from the standard Gaussian distribution, and $\tilde{X}_{0}^{\bm{a}^k, \lambda} = X_0$ is the initial state of the agent at time $t = 0$. Note that the single-realization representation of the state process coincides with an agent trial in the training procedure under the RL framework.

When solving the BSDE, we use the single-realization of the state sample $\{\tilde{X}_{n+1}^{\bm{a}^k, \lambda}\}_{n=1}^{N_T}$ generated by \eqref{SDE-sample} to represent the state process in the FBSDE system.  In this way, we rewrite the numerical scheme for the BSDE and obtain the following sample-wise approximation for the adjoint process $Y$
\begin{equation}\label{BSDE-sample}
\tilde{Y}_{n}^{\bm{a}^k, \lambda} =   \tilde{Y}_{n+1}^{\bm{a}^k, \lambda} + b_x(t_{n+1}, \tilde{X}^{\bm{a}^k, \lambda}_{n+1}, \bm{a}^{k}_{t_{n+1}}, \lambda)^{\top} \tilde{Y}^{\bm{a}^k, \lambda}_{n+1} + f^{\lambda}_x(t_{n+1}, \tilde{X}^{\bm{a}^k, \lambda}_{n+1}, \bm{a}^{k}_{t_{n+1}})^{\top} \Delta t_n, \quad n = N_T - 1, \cdots, 0,
\end{equation}
where $\tilde{Y}_{n}^{\bm{a}^k, \lambda}$ and $\tilde{Y}_{n+1}^{\bm{a}^k, \lambda}$ are approximations for $Y_{t_n}^{\bm{a}^k, \lambda}$ and $Y_{t_{n+1}}^{\bm{a}^k, \lambda}$ corresponding to the state samples $\tilde{X}_{n}^{a^k, \lambda}$ and $\tilde{X}_{n+1}^{a^k, \lambda}$, respectively, i.e. $\tilde{Y}_{n}^{\bm{a}^k, \lambda} = \tilde{Y}_{n}^{\bm{a}^k, \lambda}(\tilde{X}_{n}^{\bm{a}^k, \lambda})$ and $\tilde{Y}_{n+1}^{\bm{a}^k, \lambda} = \tilde{Y}_{n+1}^{\bm{a}^k, \lambda}(\tilde{X}_{n+1}^{\bm{a}^k, \lambda})$, and we have used stochastic approximation to approximate conditional expectations in Eq. \eqref{Scheme-Full:FBSDEs} as
$$
\tilde{Y}_{n+1}^{\bm{a}^k, \lambda}(\tilde{X}_{n+1}^{\bm{a}^k, \lambda}) \approx \E_n^{X^k}\big[ Y_{n+1}^{\bm{a}^k, \lambda}\big]
$$
and
$$
\begin{aligned}
& b_x(t_{n+1}, \tilde{X}^{\bm{a}^k, \lambda}_{n+1}, \bm{a}^{k}_{t_{n+1}}, \lambda)^{\top} \tilde{Y}^{\bm{a}^k, \lambda}_{n+1} + f^{\lambda}_x(t_{n+1}, \tilde{X}^{\bm{a}^k, \lambda}_{n+1}, \bm{a}^{k}_{t_{n+1}})^{\top} \\
\approx & \quad \E_n^{X^k} \Big[ b_x(t_{n+1}, X^{\bm{a}^k, \lambda}_{n+1}, \bm{a}^{k}_{t_{n+1}}, \lambda)^{\top} Y^{\bm{a}^k, \lambda}_{n+1} + f^{\lambda}_x(t_{n+1}, X^{\bm{a}^k, \lambda}_{n+1}, \bm{a}^{k}_{t_{n+1}})^{\top} \Big] .
\end{aligned}
$$

Then, we use sample-wise approximations introduced in \eqref{SDE-sample} and \eqref{BSDE-sample} to represent the stochastic processes $X$ and $Y$ in the gradient $\nabla J_{\bm{a}}$ and get the following sample-wise approximation for the gradient 
\begin{equation}\label{Gradient-sample}
\nabla \tilde{J}_{\bm{a}} (\bm{a}^{k}_{t_n})= b_{\bm{a}}(t_n, \tilde{X}_{n}^{\bm{a}^k, \lambda}, \bm{a}^{k}_{t_n}, \lambda)^{\top} \tilde{Y}_{n}^{\bm{a}^k, \lambda} + f^{\lambda}_{\bm{a}}(t_n, \tilde{X}_{n}^{\bm{a}^k}, \bm{a}^{k}_{t_n}, \lambda)^{\top}.
\end{equation}

As a result, the fully calculated gradient descent scheme \eqref{Approx:GD} becomes the following stochastic gradient descent scheme 
\begin{equation}\label{Approx:SGD}
\bm{a}_{t_n}^{k+1} = \bm{a}_{t_n}^{k} - \eta_k \nabla \tilde{J}_{a} (\bm{a}_{t_n}^{k}), \qquad n = 0, 1, 2, \cdots, N_T -1, \quad k = 0, 1, 2, \cdots, K-1.
\end{equation}
Although the sample-wise approximator $\{\tilde{Y}^{\bm{a}^k, \lambda}_{n}\}_{n=1}^{N_T}$ cannot provide a comprehensive representation for the adjoint process $\{Y^{\bm{a}^k, \lambda}_{t_n}\}_{n=1}^{N_T}$ in the state space since the conditional expectations are only approximated by single-realization of state samples, it's important to point out that $Y$ only appears in the gradient process under expectation, and the primary contribution of $Y$ is to incorporate the differential dynamics of the adjoint process into the gradient and guide the search of the optimal control (i.e. policy) \cite{Bao_EAJAM20}. Therefore, our sample-wise approximation for $Y$ can already embed information of the differential dynamics into the gradient process, and the rationale of applying stochastic approximation in stochastic gradient descent can be used to justify the sample-wise solver for the FBSDE system (see \cite{Bao_SNN_22} for some related analysis).

\vspace{0.5em}

In this work, we name the sample-wise numerical solver introduced by schemes \eqref{SDE-sample} - \eqref{Approx:SGD} the ``backward action learning (BAL)'' method  for solving the stochastic optimal control problem, and such a BAL framework constitutes the key mechanism of our SMP type approach for policy improvement in reinforcement learning.

\vspace{0.5em}

%It's also worthy to mention that the above BAL approach with single sample representation of the state variables is well-aligned with the RL scenario since the agent takes actions and interacts with the environment, which will create a single state trajectory that is corresponding to a simulated state sample. 

\subsection{Combine direct filter with backward action learning for solving the reinforcement learning problem}\label{BAL-RL}
Now, we combine the direct filter based exploration method with the BAL based exploitation method and construct a numerical algorithm to solve the RL problem. Since the key of the numerical recipe of our SMP approach for solving the RL problem is the BAL method for the stochastic optimal control problem, in the rest of this paper we will also call our SMP based RL solver the BAL method for convenience of presentation.

\vspace{0.5em}

To proceed, we assume that we can instantly receive the state of the agent during the training procedure. With online reception of the agent state, instead of applying the direct filter to estimate the environment parameter after each training episode to update our understanding of the environment with the information of the entire agent trial trajectory (as we introduced in Section \ref{Direct-Filter}), in the BAL algorithm for solving the RL problem we implement the direct filter dynamically at each time step in each training episode. This time-dependent parameter estimation implementation will allow us to better utilize the state information of the agent. As a result, we can more frequently update the environment information and therefore more sufficiently explore the environment.

Specifically, we let $\{\zeta_{k+1, t_0}^{(q)}\}_{q=1}^Q = \{\zeta_k^{(q)}\}_{q=1}^Q$ be the initial parameter particles at the beginning of the $k+1$-th training episode, where $\{\zeta_k^{(q)}\}_{q=1}^Q$ formulates an empirical distribution for the conditional PDF $p(\lambda_k | \mathcal{X}_k)$ of the environment parameter $\lambda$ after the $k$-th training episode. Assuming that we have a set of (equally-weighted) parameter particles $\{\zeta_{k+1, t_n}^{(q)}\}_{q=1}^Q$ at time step $t_n$ under the temporal partition $\Pi_{N_T}$, we use the following zero-dynamics scheme, which is similar \eqref{parameter:prediction}, to generate a set of predicted parameter particles  $\{\tilde{\zeta}_{k+1, t_{n+1}}^{(q)}\}_{q=1}^Q$ for time step $t_{n+1}$ in the $k+1$-th training episode
\begin{equation}\label{parameter:prediction:time}
\tilde{\zeta}_{k+1, t_{n+1}}^{(q)} =\zeta_{k+1, t_n}^{(q)} + \xi_{k+1, t_n}^{(q)}, \qquad q = 1, 2, \cdots, Q.
\end{equation}

Then, we discretize the state equation at time step $t_n$ corresponding to the parameter particles as follows
\begin{equation}\label{Agent:State}
\tilde{X}_{n+1}^{\bm{a}^k, \tilde{\zeta}_{k+1, t_{n+1}}^{(q)}} = \tilde{X}_{n}^{\bm{a}^k, \tilde{\zeta}_{k+1, t_{n}}^{(q)}} + b(t_n, \tilde{X}^{\bm{a}^k, \tilde{\zeta}_{k+1, t_{n}}^{(q)}}_{n}, \bm{a}_{t_n}^k, \tilde{\zeta}_{k+1, t_{n+1}}^{(q)}) \Delta t_n + \sigma_{t_n} \sqrt{\Delta t_n} \epsilon_n^{(q)}, \quad q = 1, 2, \cdots, Q.
\end{equation}
Note that the Gaussian random samples $\{\xi_{k+1, t_n}^{(q)}\}_{q=1}^{Q}$ add artificial noise to the parameter cloud $\{\tilde{\zeta}_{k+1, t_{n+1}}^{(q)}\}_{q=1}^Q$, and it can provide a natural mechanism that encourages the agent to explore the environment. Then, we use $X^{\bm{a}^k, \lambda}$ to denote the trajectory of the $k+1$-th agent trial, which interacts with the real environment, and we compare the simulated state-parameter samples $\tilde{X}_{n+1}^{\bm{a}^k, \tilde{\zeta}_{k+1, t_{n+1}}^{(q)}}$ (generated by Eq. \eqref{Agent:State}) with the real agent state $X_{t_{n+1}}^{\bm{a}^k, \lambda}$ to derive a likelihood value for each predicted parameter particle $\tilde{\zeta}_{k+1, t_{n+1}}^{(q)}$, i.e.
\begin{equation}\label{Likelihood:particles}
\begin{aligned}
& \qquad p(X_{t_{n+1}}^{\bm{a}^k, \lambda}, f^{\lambda}({t_{n+1}}, X_{t_{n+1}}^{\bm{a}^k, \lambda}, \bm{a}^k) \big| \tilde{\zeta}_{k+1, t_{n+1}}^{(q)}) \\ 
\sim & \ \exp\Big( - \f{\big( X_{t_{n+1}}^{\bm{a}^k, \lambda} -  X_{n+1}^{\bm{a}^k,  \tilde{\zeta}_{k+1, t_{n+1}}^{(q)}}\big)^2 + \big( f^{\lambda}({t_{n+1}}, X_{t_{n+1}}^{\bm{a}^k, \lambda}, \bm{a}^k)-  f^{\tilde{\zeta}_{k+1, t_{n+1}}^{(q)}}({t_{n+1}}, X_{n+1}^{\bm{a}^k,  \tilde{\zeta}_{k+1, t_{n+1}}^{(q)}}, \bm{a}^k)\big)^2}{2 \big(\delta_k^2 + (\sigma_{t_n} \Delta t)^2 \big)}\Big).
\end{aligned}
\end{equation}
As a result, the particle-weight pairs $\{\tilde{\zeta}_{k+1, t_{n+1}}^{(q)}, \omega_{k+1, t_{n+1}}^{(q)}\}_{q=1}^Q$, where the weight value is assigned as $\omega_{k+1, t_{n+1}}^{(q)} = p(X_{t_{n+1}}^{\bm{a}^k, \lambda}, , f^{\lambda}({t_{n+1}}, X_{t_{n+1}}^{\bm{a}^k, \lambda}, \bm{a}^k) \big| \tilde{\zeta}_{k+1, t_{n+1}}^{(q)})/C$ with an appropriate normalization factor $C$, form a weighted empirical distribution for the posterior distribution of the parameter at time instant $t_{n+1}$ in the $k+1$-th training episode.  To avoid the degeneracy issue in each training episode, we resample all the particles based on the weighted pairs $\{\tilde{\zeta}_{k+1, t_{n+1}}^{(q)}, \omega_{k+1, t_{n+1}}^{(q)}\}_{q=1}^Q$ and generate a set of equally-weighted particles $\{\zeta_{k+1, t_{n+1}}^{(q)}\}_{q=1}^Q$ for the next time step. 

In this way, at the terminal time $T$ we have carried out $N_T-1$ times parameter estimation procedure, which can effectively incorporate the information of the agent state at each time instant in the $k+1$-th training episode, and the particle set $\{\zeta_{k+1}^{(q)}\}_{q=1}^Q := \{\zeta_{k+1, t_{N_T}}^{(q)}\}_{q=1}^Q$ provides our ``best'' understanding of the environment after considering the state of the $k+1$-th agent trial. Then, we use the mean of the parameter particles $\{\zeta_{k+1}^{(q)}\}_{q=1}^Q$, denoted by $\bar{\lambda}_{k+1}$, as our estimate for the environment parameter in the $k+1$-th training episode. 
%\begin{equation}\label{Estimate:lambda}
%\bar{\lambda}_{k+1} := \f{1}{Q}\sum_{q=1}^Q \zeta_{k+1}^{(q)}
%\end{equation}

With the updated estimate $\bar{\lambda}_{k+1}$ for the environment parameter, we follow the BAL algorithm discussed in Section \ref{BAL-SOC} to improve the policy in the exploitation procedure. Recall that the main theme of the BAL method is to use a single-realization of the state trajectory as a stochastic approximation to the state process when approximating expectations. Since the real agent trajectory $X^{\bm{a}^k, \lambda}$ that follows policy $\bm{a}^k$ already provides a path of the agent state, we can use the real agent trial state at the temporal partition points, i.e. $\{X_{t_n}^{\bm{a}^k, \lambda}\}_{n=1}^{N_T}$, to replace the single-realization simulated state trajectory $\{\tilde{X}_{n}^{\bm{a}^k, \lambda}\}_{n=1}^{N_T}$ (introduced in Eq. \eqref{SDE-sample}) for the BAL optimal control solver. On the other hand, the adjoint equation, i.e., the BSDE, is still driven by the current estimated parameter $\bar{\lambda}_{k+1}$ except that the forward process is replaced by the real agent state, and we introduce the following sample-wise solution for the adjoint BSDE in the $k+1$-th training episode:
\begin{equation}\label{BSDE-sample:RL}
\begin{aligned}
\bar{Y}_{n}^{\bm{a}^k, \bar{\lambda}_{k+1}} =  & \bar{Y}_{n+1}^{\bm{a}^k, \bar{\lambda}_{k+1}} + b_x(t_{n+1}, X^{\bm{a}^k, \lambda}_{t_{n+1}}, \bm{a}^{k}_{t_{n+1}}, \bar{\lambda}_{k+1})^{\top} \bar{Y}^{\bm{a}^k, \bar{\lambda}_{k+1}}_{n+1} \\
& \quad + f^{\bar{\lambda}_{k+1}}_x(t_{n+1}, X^{\bm{a}^k, \lambda}_{t_{n+1}}, \bm{a}^{k}_{t_{n+1}})^{\top} \Delta t_n, \qquad n = N_T - 1, \cdots, 0.
\end{aligned}
\end{equation}
Note that we have used $\bar{\lambda}_{k+1}$ as our estimated environment parameter in $b$ and $f$ when $\lambda$ is explicitly needed in the numerical scheme for the BSDE, and the single-realization representation of the forward SDE is given by the agent trial trajectory $\{X_{t_n}^{\bm{a}^k, \lambda}\}_{n=1}^{N_T}$.  

Then, we derive the sample-wise approximation for the gradient as follows
\begin{equation}\label{Gradient-sample:RL}
\nabla \bar{J}_{\bm{a}} (\bm{a}^{k}_{t_n})= b_{\bm{a}}(t_n, X^{\bm{a}^k, \lambda}_{t_{n}}, \bm{a}^{k}_{t_n}, \bar{\lambda}_{k+1})^{\top} \bar{Y}_{n}^{\bm{a}^k, \bar{\lambda}_{k+1}} + f^{\bar{\lambda}_{k+1}}_{\bm{a}}(t_n, X^{\bm{a}^k, \lambda}_{t_{n}}, \bm{a}^{k}_{t_n})^{\top},
\end{equation}
and we carry out the following stochastic gradient descent iteration to improve the policy
\begin{equation}\label{Approx:SGD:RL}
\bar{\bm{a}}_{t_n}^{k+1} = \bar{\bm{a}}_{t_n}^{k} - \eta_k \nabla \bar{J}_{\bm{a}} (\bm{a}_{t_n}^{k}), \qquad n = 0, 1, 2, \cdots, N_T -1.
\end{equation}

\vspace{0.5em}
Different from the classic stochastic optimal control problem, which aims to find the optimal control with an explicitly given environment, the agent in the RL problem needs to explore the environment corresponding to the state space. In the above BAL approach for solving the RL problem, the exploration is implemented through the parameter estimation procedure, and the artificial noise added to the pseudo parameter dynamics (as described in Eq. \eqref{parameter:prediction:time}) allows the agent to explore. In order to explore more actively and detect other possibilities, we adopt the $\epsilon$-greedy type exploration algorithm and provide a mechanism for the agent to randomly explore the environment. 

Specifically, for the suboptimal policy $\bm{a}^k$ that we obtained in the $k$-th training episode,  we let the updated policy for the $k+1$-th training episode as follows
\begin{equation}\label{epsilon-greedy}
\bm{a}^{k+1}_{t_n} = \left\{
\begin{aligned}
&\bar{\bm{a}}^{k+1}_{t_n}, \quad &\text{with probability $1-\epsilon$}\\
&\bar{\bm{a}}^{k+1}_{t_n} + \theta_n &\text{with probability $\epsilon$}\\
\end{aligned} \right., \qquad n = 0, 1, 2, \cdots, N_T - 1,
\end{equation}
where $\theta_n \sim N(0, \Sigma)$ is a user defined noise level that brings random perturbations to the policy with the size of a pre-determined covariance $\Sigma$, and $0 \leq \epsilon < 1$ is the probability of implementing the enhanced policy exploration.

\hspace{0.5em}

It's important to point out that the estimated optimal policy $\bm{a}^K$ obtained through the above explorative stochastic gradient descent optimization procedure only gives a deterministic policy process, which is based on the initial state $X_0$. In this way, only the policy at time $t_0 = 0$, i.e. $\bm{a}_{t_0}^K = \bm{a}_{t_0}^{K}(X_0)$, is a state-dependent action that reflects the actual state-to-action map, and the policy process beyond time $t_0$ only gives an estimate based on the expected future behavior of the agent.
In order to let the agent take optimal actions based on its current state beyond the initial time $t_0$ as desired in the RL problem, we need to screen the state space and calculate the optimal policy corresponding to the state at any time instant $t_n$.  In what follows, we shall modify the BAL method to generate a time/state-dependent optimal policy process.

Note that the state of the agent at time $t_{n+1}$ only depends on the previous states and actions. Therefore, we don't need to search the optimal policy in the entire state space.  At the initial time $t_0$,  assume that we have calculated the policy process $\{\bm{a}^K_{t_n}\}_{n=0}^{N_T - 1}$ through schemes \eqref{BSDE-sample:RL} - \eqref{epsilon-greedy} based on the estimated environment parameter obtained by the direct filter method. With the initial state $X_0$ and the optimal action at time $t_0$, i.e. $\bm{a}_{t_0}^K = \bm{a}_{t_0}^K(X_{t_0})$, we carry out the following Euler-Maruyama scheme to generate $M$ samples for $X_{1}^{\bm{a}^K, \lambda}$
\begin{equation}\label{Predict:region}
\tilde{X}_{1}^{\bm{a}^k, \lambda, (m)} = X_0 + b(t_n, X_0, \bm{a}_{t_0}^k(X_0), \lambda) \Delta t_0 + \sigma_{t_0} \sqrt{\Delta t_0} \epsilon_0^{(m)}, \quad m = 1, 2, \cdots, M,
\end{equation}
where $\{\epsilon_0^{(m)}\}_{m=1}^M$ are $M$ samples drawn from the standard Gaussian distribution, and the state samples $\{\tilde{X}_{1}^{\bm{a}^k, \lambda, (m)}\}_{m=1}^M$ characterize the state variable $X_{t_1}^{\bm{a}^K, \lambda}$ in the state space. Then, we let $\mathcal{D}_{t_1}$ be the region in the state space that covers all the state samples $\{\tilde{X}_{1}^{\bm{a}^k, \lambda, (m)}\}_{m=1}^M$. Apparently, if we run Eq.\eqref{Predict:region} with a large enough number $M$, the simulated samples $\{\tilde{X}_{1}^{\bm{a}^k, \lambda, (m)}\}_{m=1}^M$ would provide very reliable predictions for the future state $X_{t_1}^{\bm{a}^K, \lambda}$, and the corresponding region $\mathcal{D}_{t_1}$ would be large enough to cover the agent state at time $t_1$ if the agent takes the optimal action at time $t_0$. In this way, actions corresponding to state points in the region $\mathcal{D}_{t_1}$ at time $t_1$ are needed. To proceed, we introduce a set of state points, denoted by $\mathcal{X}_{t_1}$, as a spatial discretization for the state region $\mathcal{D}_{t_1}$. Then, we carry out the same BAL algorithm \eqref{BSDE-sample:RL} - \eqref{epsilon-greedy} from initial time $t_1$ to terminal time $t_{N_T}$, and the initial state is chosen among the state points in $\mathcal{X}_{t_1}$, i.e. $X_{t_1}^{\bm{a}^K, \lambda}  = x\in \mathcal{X}_{t_1}$. As a result, we obtain a set of state-dependent optimal actions $\{\bm{a}^{K}_{t_1}(x)\}_{x\in \mathcal{X}_{t_1}}$ at time instant $t_1$, and we use  $\{\bm{a}^{K}_{t_1}(x)\}_{x\in \mathcal{X}_{t_1}}$ as our policy variable (or the policy table) to guide the agent at time $t_1$. Similarly, following the above procedure, if we start from state points in $\mathcal{X}_{t_n}$ with their optimal actions  $\{\bm{a}^{K}_{t_n}(x)\}_{x \in \mathcal{X}_{t_n}}$, we can determine the state region $\mathcal{D}_{t_{n+1}}$ and the state points $\mathcal{X}_{t_{n+1}}$ at time $t_{n+1}$. Then, we can compute the optimal policy $\{\bm{a}^{K}_{t_{n+1}}(x)\}_{x \in \mathcal{X}_{t_{n+1}}}$ corresponding to the state points in $\mathcal{X}_{t_{n+1}}$ by using the BAL method. As a result, we can adaptively calculate the optimal policy over the state space along the temporal partition $\Pi_{N_T}$.

We also want to point out that although the training procedure for the optimal policy needs to be repeated carried out step-by-step in time, the direct filter based exploration can provide good understanding of the environment through the exploration procedure even at the first time instant $t_0$ since the agent trial trajectories already explored the environment with a well-designed deterministic policy. Although the uncertainty in the state model may lead the agent to different possible paths in the future, the deterministic optimal policy can be close to the stochastic optimal policy. The training procedure for the optimal policy after time instant $t_0$ mainly incorporates stochasticity to the policy process so that the agent can take appropriate actions corresponding to the random state due to the uncertainty generated by the Brownian motion $W$ in the state dynamics. In other words, we can get better and more complete estimate for the environment parameter as time. This makes the RL problem that we consider in this work become more and more like a classic stochastic optimal control problem except that we allow the agent to test the BAL designed policy as trials.

\subsection{Summary of the algorithm}
To summarize our BAL algorithm for solving the RL problem, we provide a pseudo algorithm in Table \ref{Algorithm}.
\vspace{1em}

\begin{table}\caption{Summary of the algorithm}\label{Algorithm}
%\vspace{1em}
\centering
\begin{tabular} {p{0.9\textwidth}}
\hline\noalign{\smallskip}
{\bf Algorithm}: {\em Backward action learning method for reinforcement learning}\\
\noalign
{\smallskip}\hline
\noalign{\smallskip}
\vspace{-0.1cm}
\begin{spacing}{1.1}
\begin{algorithmic}
%\item[Initialize] the number of particles $Q$ for the direct filter, the learning rate $\{\eta_{k}\}_{k=1}^K$ and the total number of training episodes $K \in \mathbb{N}$. Generate parameter particles $\{\zeta_1^{(q)}\}_{q=1}^Q$ for the initial guess of the environment parameter and introduce an initial policy $\bm{a}^{1}$.
\item[Set up] the initial state $X_0$, environment guess particles $\{\zeta_{0}^{(q)}\}_{q=1}^Q$, initial policy $\bf{a}^0$, and choose the user defined constants: $K$, $Q$, $M$, $\{\eta_k\}_{k=1}^K$, and $\epsilon$.
\vspace{0.3em}
\item[\textbf{for}] $n = 0, 1, 2, \cdots, N_T - 1$
\begin{itemize}
\item[Let ] $X_{t_n}^{\bm{a}^0} = x \in \mathcal{X}_{t_n}$ be the initial state at time $t_n$, and obtain the optimal policy $\bm{a}^{K}_{t_n} (x)$ for each $x \in \mathcal{X}_{t_n}$ as follows:

\hspace{-3em} \textbf{while} $k =0, 1, 2, \cdots, K$, \textbf{do} \\
\begin{description}
\vspace{-0.75em}
\item[$\bullet$] \textbf{Implement policy $\bm{a}_{t_n:N_T-1}^k$ and obtain the agent trial state trajectory $X_{t_n:N_T}^{\bm{a}^k, \lambda}$.}
\vspace{-1.25em}
\item[$\bullet$] \textbf{Implement the direct filter method to explore:}

\item[$\Box$] Let $\{\zeta_{k+1, t_{n}}^{(q)}\}_{q=1}^Q = \{\zeta_{k}^{(q)}\}_{q=1}^Q$ be the particles for the estimated environment parameter

\textbf{for} $l = n, n+1, \cdots, N_T - 1$
\begin{description}
\vspace{-0.75em}
\item[-] Generate a set of $Q$ predicted parameter particles $\{\tilde{\zeta}_{k+1, t_{l+1}}^{(q)}\}_{q=1}^Q$ through the pseudo parameter dynamics Eq. \eqref{parameter:prediction:time};
\vspace{-0.75em}
\item[-]  Generate $Q$ versions of the state samples $\{\tilde{X}_{l+1}^{\bm{a}^k, \tilde{\zeta}_{k+1, t_{l+1}}^{(q)}}\}_{q=1}^Q$ from the approximation scheme Eq. \eqref{Agent:State} for the next state stage;
\vspace{-0.75em}
\item[-] Calculate likelihood values for $\{\tilde{\zeta}_{k+1, t_{l+1}}^{(q)}\}_{q=1}^Q$ as particle weights $\{ \omega_{k+1, t_{l+1}}^{(q)}\}_{q=1}^Q$ by comparing $\{\tilde{X}_{l+1}^{\bm{a}^k, \tilde{\zeta}_{k+1, t_{l+1}}^{(q)}}\}_{q=1}^Q$ with the real agent state $X_{t_{l+1}}^{\bm{a}^k, \lambda}$ through Eq. \eqref{Likelihood:particles};
\vspace{-0.75em}
\item[-] Resample the particle-weight pairs $\{\tilde{\zeta}_{k+1, t_{l+1}}^{(q)}, \omega_{k+1, t_{l+1}}^{(q)}\}_{q=1}^Q$ to generate a set of equally weighted particles $\{\zeta_{k+1, t_{l+1}}^{(q)}\}_{q=1}^Q$.
\vspace{-0.75em}
\end{description}
\textbf{end for}

\item[$\Box$]  Set  $\{\zeta_{k+1}^{(q)}\}_{q=1}^Q = \{\zeta_{k+1, t_{N_T}}^{(q)}\}_{q=1}^Q$ to initialize the next training episode, and let $\bar{\lambda}_{k+1} = \f{1}{Q}\sum_{q=1}^Q \zeta_{k+1}^{(q)}$ be the estimated environment parameter in the training episode $k+1$.
\item[$\bullet$] \textbf{Implement the BAL method to exploit:}
\begin{description}
\vspace{-0.75em}
\item[-] Solve the adjoint BSDE through the numerical BSDE scheme Eq. \eqref{BSDE-sample:RL} by using the agent state trajectory $X_{t_n:N_T}^{\bm{a}^k, \lambda}$ and the estimated environment parameter $\bar{\lambda}_{k+1}$;
\vspace{-0.75em}
\item[-] Calculate the sample-wise approximator $\nabla \bar{J}_{\bm{a}}$ for the gradient introduced in Eq. \eqref{Gradient-sample:RL};
\vspace{-0.75em}
\item[-] Carry out stochastic gradient descent scheme \eqref{Approx:SGD:RL}, to get the improved policy $\bar{\bm{a}}^{k+1}$.
\vspace{-1.5em}
\end{description}
\item[$\bullet$] Use $\epsilon$-greedy method \eqref{epsilon-greedy} to enhance exploration and obtain $\bm{a}^{k+1}$.
\end{description}
\hspace{-3em}\textbf{end while}
\end{itemize}
\vspace{-0.5em}
Generate the next state region $\mathcal{D}_{t_{n+1}}$ from the state points $\mathcal{X}_{t_n}$ and the policy table $\{\bm{a}^{K}(x)\}_{x \in \mathcal{X}_{t_n}}$ at time $t_n$ through scheme \eqref{Predict:region}. Discretize $\mathcal{D}_{t_{n+1}}$ and create state points $\mathcal{X}_{t_{n+1}}$.
\item[\textbf{end for}]
\end{algorithmic}
\vspace{-1.2em}
\end{spacing}\\
\hline
\end{tabular}
\end{table}

%\newpage

\section{Numerical experiments}\label{Numerics}
In this section, we use three numerical examples to demonstrate the performance of our BAL method for solving the RL problem. 

\subsection{Example1: Classic linear-quadratic control with a hidden environment parameter.}

In the first example, we solve a classic linear-quadratic stochastic optimal control problem, in which the state model contains a hidden parameter that represents the unknown in the environment. The main purpose of demonstrating a linear-quadratic control example is that the optimal control can be explicitly derived, hence we can use this example to present the accuracy of our BAL method by comparing with the analytically derived solution. On the other hand, the unknown parameter in the state model requires an exploration procedure to determine the environment, which also makes the control problem in this example an RL problem.

\vspace{0.5em}

To proceed, we consider an agent, whose state is formulated by the following  2-dimensional linear stochastic dynamical system:
\begin{equation}\label{State:LQ}
dX^{\bm{a}, \lambda}_t = (A(t) X^{\bm{a}, \lambda}_t + B \bm{a}_t) dt + \sigma dW_t,
\end{equation}
where $X^{\bm{a}, \lambda}_t \in \mathbb{R}^2$ is the state of the agent;  $A(t) = [\lambda \sin t, 0; 0, \cos t]$ is the drift coefficient for state $X$, which contains an unknown parameter $\lambda$, and we choose $\lambda = 2$ in our numerical experiments in this example; $B = (0.5, 0.5)^{\top}$ is a 2-dimensional constant vector as the coefficient of the scalar policy term $\bm{a}_t$; $\sigma = [0.1, 0; 0, 0.1]$ is the diffusion coefficient for the stochastic integral driven by the 2-dimensional Brownian motion $W$.  The cost functional, which is equivalent to the penalty in the RL problem, is in the following quadratic form
$$J(\bm{a}) = E \Big[\frac{1}{2}\int_{0}^{T} \Big(\lan Q X^{\bm{a}, \lambda}_t, X^{\bm{a}, \lambda}_t\ran+\lan R\bm{a}_t, \bm{a}_t \ran \Big) dt + \frac{1}{2}\lan FX_T, X_T\ran\Big],$$
where $Q = I_2$ and $F = I_2$ are given constant matrices, and we let $R = 1$.  In the numerical experiments, we introduce a temporal partition over $[0, 1]$ with the uniform step-size $\Delta t = 0.01$, i.e. $N_T = 100$. 

We first show the performance of our direct filter parameter estimation method for the purpose of exploration in the RL task, and we use $Q = 100$ particles to describe the empirical distribution for the unknown parameter $\lambda$.  The initial state of the agent is chosen as $X^{\bm{a}, \lambda}_0 = (6, -2)^{\top}$, and we assume that the initial guess for the environment parameter is $-2$.
%\vspace{-0.3em}
\begin{figure}[h!]
\begin{center}
\includegraphics[scale = 0.9]{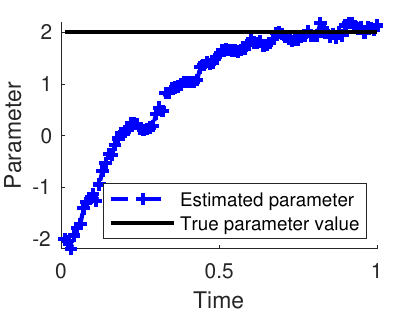}
%\vspace{-0.25em}
\end{center}
%\vspace{-0.75em}
\caption{Example 1. Parameter estimation with respect to time in the first episode. }\label{LQ_A_est} %\vspace{-0.5em}
\end{figure}
In Figure \ref{LQ_A_est}, we present the estimated parameter with respect to time in the first training episode, where the solid black line shows the true environment parameter value $\lambda = 2$ while the blue dashed curve marked by plus signs gives our estimated parameter values corresponding to time in the training episode, i.e. $k=1$.  We can see from this figure that the direct filter method can quickly capture the true parameter even in the first training episode in this 1-dimensional state-independent parameter estimation case. In the next two examples (Example 2 and Example 3), we will consider more complicated situations with state-dependent environment parameters, which can be challenging for standard reinforcement learning techniques.

\begin{figure}[h!] %\vspace{-0.75em}
\begin{center}
\subfloat[$50$ episodes]{\includegraphics[scale = 0.55]{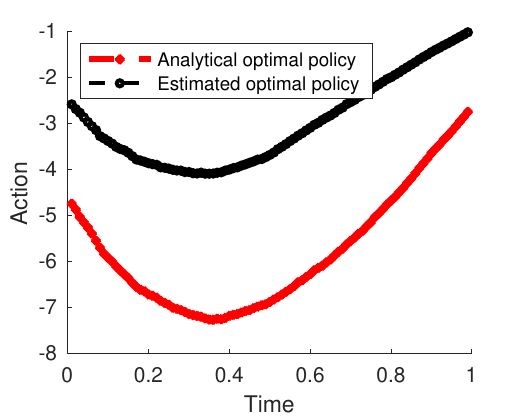} } \qquad
\subfloat[$100$ episodes]{\includegraphics[scale = 0.55]{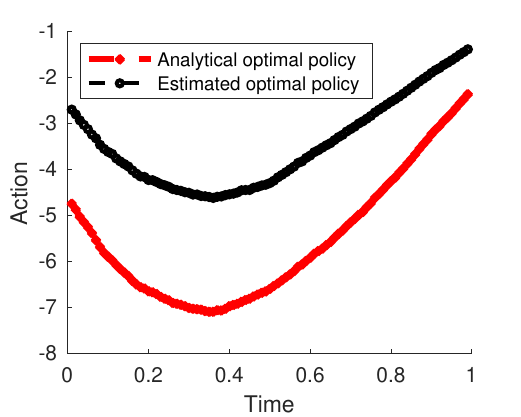}}\\
\subfloat[$500$ episodes]{\includegraphics[scale = 0.55]{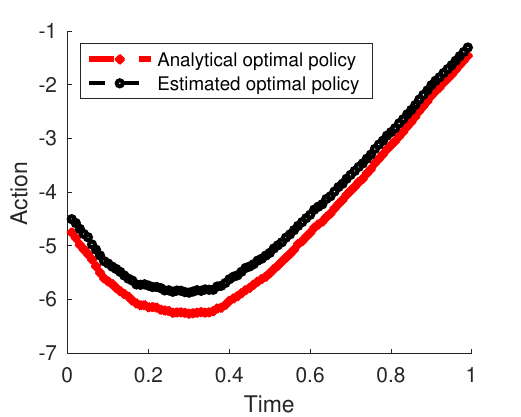}} \qquad
\subfloat[$1,000$ episodes]{\includegraphics[scale = 0.55]{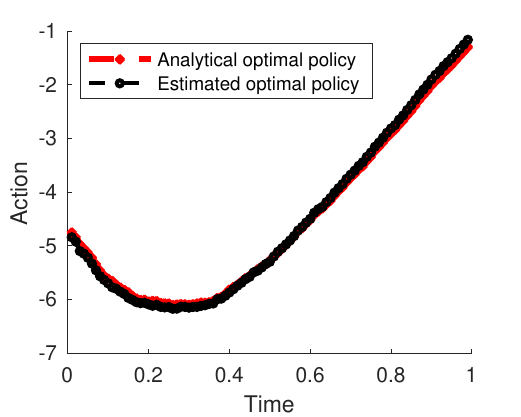}}
\end{center}
%\vspace{-0.75em}
\caption{Example 1. Comparison between the estimated optimal policy (actions) and the analytical optimal control. }\label{LQ_iter}
\end{figure}
With the accurately estimated environment, we present the performance of policy estimation with respect to the number of training episodes in Figure \ref{LQ_iter}. In subplots (a), (b), (c), and (d) we compare the BAL method estimated optimal policy (the black dashed curves marked by circles) with the analytical optimal policy (the black dashed curves), which is given by
\begin{equation}\label{Ex1:control}
\bm{a}^{\ast}_t = - R^{-1} B^TP(t)X^{\bm{a}, \lambda}_t,
\end{equation}
where $X_t$ is the agent state, and $P(t)$ is the unique solution of the following so-called Riccati equation corresponding to the state equation \eqref{State:LQ}
$$\frac{d P(t)}{dt} = - P(t) A(t) - A^T(t) P(t) + P(t) B R^{-1}  B^T P(t) - Q , \quad P(T) = F.$$ 
From this figure, we can see that as we increase the number of training episodes, i.e. $K =50$, $K=100$, $K=500$, and $K=1,000$ in (a), (b), (c), and (d), respectively, we obtain better and better policy estimation results. When carrying out $K=1,000$ training episodes, the estimated optimal policy is almost perfectly aligned with the analytical optimal policy. 

\begin{figure}[h!] %\vspace{-0.75em}
\begin{center}
\subfloat[$X_0 = (4, 4)^{\top}$ ]{\includegraphics[scale = 0.6]{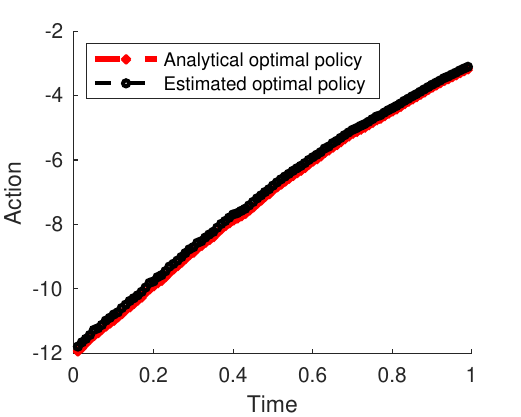} } 
\subfloat[$X_0 = (4, 1)^{\top} $]{\includegraphics[scale = 0.6]{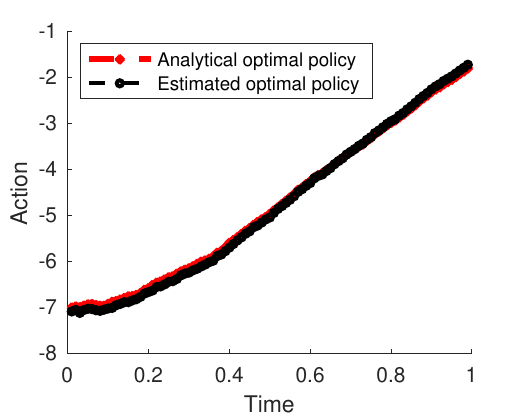}}
\subfloat[$X_0 = (6, -2)^{\top}$]{\includegraphics[scale = 0.6]{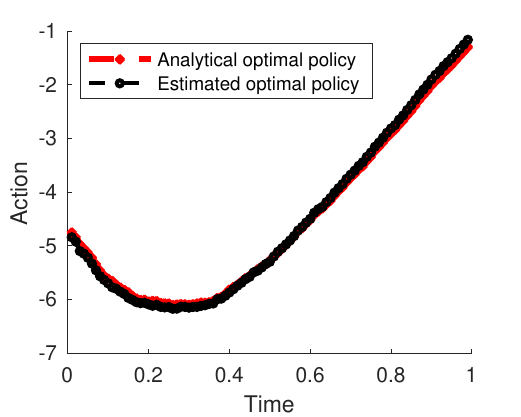}}
\end{center}
%\vspace{-0.75em}
\caption{Example 1. Policy accuracy with different initial state. }\label{LQ_state}
\end{figure}
Finally, in Figure \ref{LQ_state} we show the accuracy of our policy estimation with different initial states $X^{\bm{a}, \lambda}_0 = (4, 5)^{\top}$, $X^{\bm{a}, \lambda}_0 = (4, 1)^{\top} $, and $X^{\bm{a}, \lambda}_0 = (6, -2)^{\top}$. We can see that the BAL method constantly provides accurate policy estimation results.

\subsection{Example 2: Reinforcement learning for atomic level manufacture.}

In the second example, we solve a mathematically modified material science problem that motivated us to develop such an SMP based method to solve the RL problem with parameterized environment. We want to use this problem as an example to show that there are application problems which require physics knowledge to be incorporated into the RL model, and there are needs to parameterize the environment in scientific machine learning practices.

The scientific background of the RL problem that we consider in this example is known as the ``atomic forge'', which is a technique to control and design materials at the nano scale  \cite{Kalinin-Atom}. A new nano-phase fabrication approach has been developed to utilize a scanning transmission electron microscope (STEM) to assemble and manipulate matter atom-by-atom \cite{Bao_Atom20, Ali-Atom, Bao_Carbon_22}. 
Although many practical challenges still need to be addressed to achieve the atomic forge technique from the physics aspect, in this work we focus on the mathematical problem about how to automatically control atoms and formulate an RL method to design an effective policy to move a target atom to a pre-designated location on a 2-dimensional material surface, where two-atom potential models can be applied  \footnote{Moving atoms in the 3-dimensional space would be similar from the mathematical aspect. However, it could be more challenging in physics. }.

The motion of a target atom, which is the agent in the RL problem, is mainly driven by atomic forces derived from intermolecular potentials. One of the most important intermolecular potentials is the Lennard-Jones (LJ) potential, which models soft repulsive and attractive interactions between two atoms.  In this work, we consider the following AB form of the LJ potential
$$V^{\lambda}_{LJ}(r) := \lambda \Big( \f{A}{r^{12}} - \f{B}{r^6}\Big),$$ 
where $r$ is the distance between two interacting atoms, $\lambda$ is the depth of the potential-well (usually referred to as ``dispersion energy''), and $A$, $B$ are constant values referred to as ``size of the atom''. In this work, we let $A = B = 0.5$ be pre-chosen constants, and the depth of the potential-well $\lambda$ is an unknown value to be determined. Since we try to move a target atom on a material surface, the value of $\lambda$ may vary depending on the type of the fixed background atom, which is interacting with the moving target atom in the environment.

\begin{figure}[h!]
\begin{center}
\includegraphics[scale = 0.8]{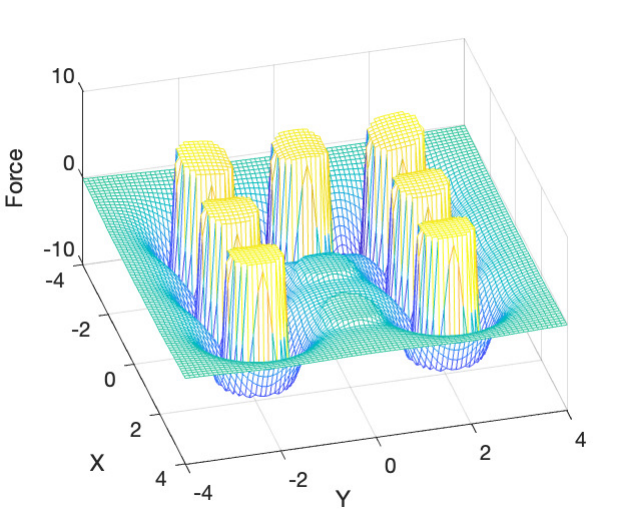}
%\vspace{-0.25em}
\end{center}
%\vspace{-0.75em}
\caption{Example 2. The atomic force on the material surface. }\label{Atomic_Force} %\vspace{-0.5em}
\end{figure}

The atomic force between atoms is in the form of the gradient of the potential $\nabla V^{\lambda}_{LJ}$. In Figure \ref{Atomic_Force}, we present a demonstration for the field of the atomic force generated by the LJ potential $V^{\lambda}_{LJ}$ corresponding to two types of altogether $9$ background atoms with depth parameters $\lambda_{deep} = 30$ and $\lambda_{shallow} = 1$. 
We can see from the figure that corresponding to the deep potential parameter, i.e. $\lambda_{deep} = 30$, both the repulsive and attractive forces are large, while on the other hand the intermolecular force is generally much small near the atoms with the shallow potential parameter, i.e. $\lambda_{shallow} = 1$. Also, the moving target atom cannot get too close to the fixed background atom due to the exponentially increased repulsive force, and the target atom can be trapped by one of the background atoms with deep potential well due to the large attractive force applied to the target atom.

\vspace{0.5em}

With the assumed intermolecular potential $V_{LJ}$ and a background atomic structure as the environment, the trajectory of the target atom can be formulated as
\begin{equation}\label{Atom-state}
dX^{\bm{a}, \lambda}_t = (-\nabla V^{\lambda}_{LJ}(r) + \bm{a}_t)dt + \sigma_t dW_t,
\end{equation}
where $\nabla V^{\lambda}_{LJ}$ is the gradient of the LJ potential $V^{\lambda}_{LJ}$, which is determined by the depth parameter $\lambda$ of the potential-well, and $r= \|X^{\bm{a}, \lambda}_t - \text{Atom}_{background}\|_2$ is the  distance between the target atom and it's closest background atom, i.e. $\text{Atom}_{background}$. The policy $\bm{a}_t$ is the control actions that we apply to drive the target atom and guide it to the pre-designated location. More specifically,  we let $\bm{a}_t : = (f_t \cos \theta_t, f_t \sin \theta_t)^{\top}$, where $f_t$ is the amount of external force that we apply to counter-effect the background atomic force caused by the LJ potentials, and $\theta_t$ is the steering angle that determines the direction of the external force.

\begin{figure}[h!]
\begin{center}
\includegraphics[scale = 0.7]{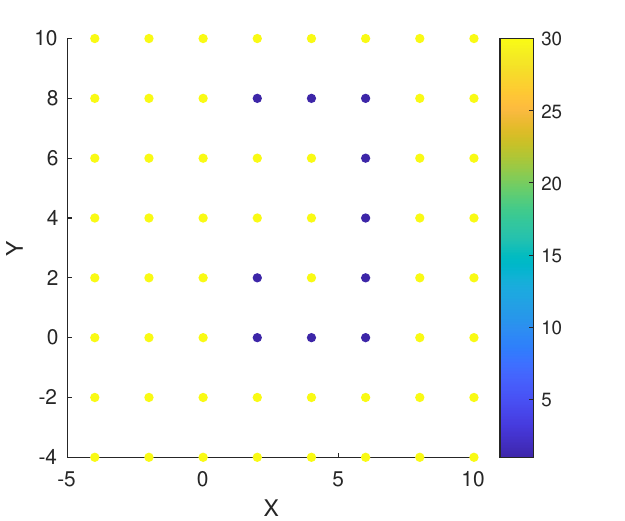}
%\vspace{-0.25em}
\end{center}
%\vspace{-0.75em}
\caption{Example 2. The true potential parameters on the material surface. }\label{Atomic_True_Parameter} %\vspace{-0.5em}
\end{figure}

To demonstrate how an RL algorithm can be applied to the atomic forge technique, we design a background atomic structure in Figure \ref{Atomic_True_Parameter}, where the blue dots show the locations of the background atoms that can generate shallow potential wells, the yellow dots show the locations of the background atoms that can generate deep potential wells, and the color bar on the right hand side shows the mapping from color to LJ potential values. 
In this example, we consider the moving target atom as the agent in the RL problem, and we choose a pre-designated destination $X_{destination} = (2, 8)^{\top}$ to be arrived at the terminal time $T = 10$. The initial state of the agent is chosen as $X_0 = (4.75, 0.75)^{\top}$, which is near a background atom that generates shallow potential well. The cost for the RL problem is defined by
\begin{equation}\label{Atom-cost}
J(\bm{a}) = \E\left[ \int_0^T |f_t| dt + F \|X^{\bm{a}, \lambda}_T - X_{destination}\|_2\right],
\end{equation}
where $F = 50$ is the amount of penalty for not being able to arrive at the destination at the terminal time $T$, and $\int_0^T |f_t| dt$ is the running cost that measures how much effort (or energy) that we make to move the agent. The RL problem is to find an optimal policy $\bm{a}^{\ast}$ that minimizes the cost $J(\bm{a})$.  It's worthy to point out that such an RL problem is quite challenging since once the target atom is ``captured'' by one of the deep potential-well atoms, it needs very large force to drive it out of the potential-well.

We first use the standard Q-Learning method with $\epsilon$-greed exploration (see \cite{Sutton_RL_2nd, Q_learning_92}) to solve the RL problem described by Eq. \eqref{Atom-state}-\eqref{Atom-cost}. To discretize the original continuous problem, we introduce a temporal partition with step-size $\Delta t = 0.2$, i.e. $N_T = 50$, and we introduce a uniform spatial partition with step-size $\Delta x = 0.1$ in the state space to generate the q-table. The policy approximation for the q-table  is chose as $\Delta f = 0.2$ for discretizing the external force $f$ and $\Delta \theta = \f{\pi}{8}$ for discretizing the steering direction $\theta$.

\begin{figure}[h!] %\vspace{-0.75em}
\begin{center}
\subfloat[$10^3$ training episodes in Q-Learning]{\includegraphics[scale = 0.55]{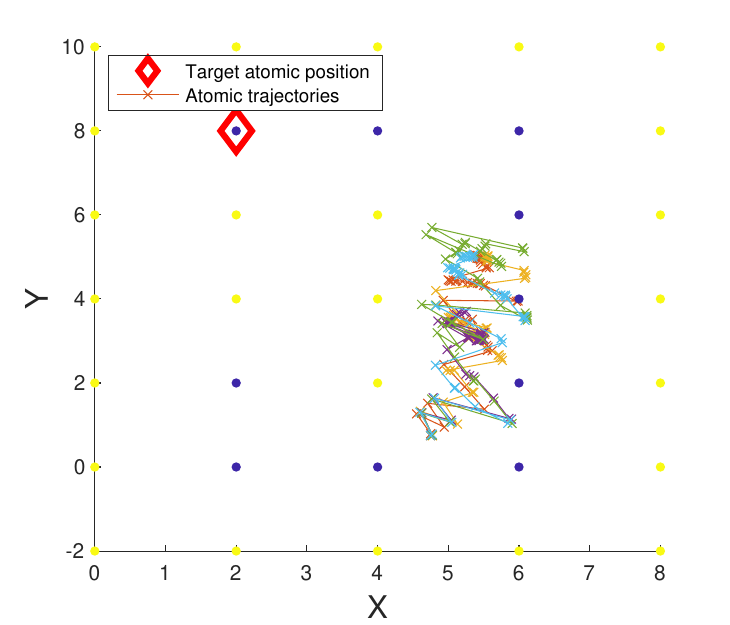} } 
\subfloat[$10^4$ training episodes in Q-Learning]{\includegraphics[scale = 0.55]{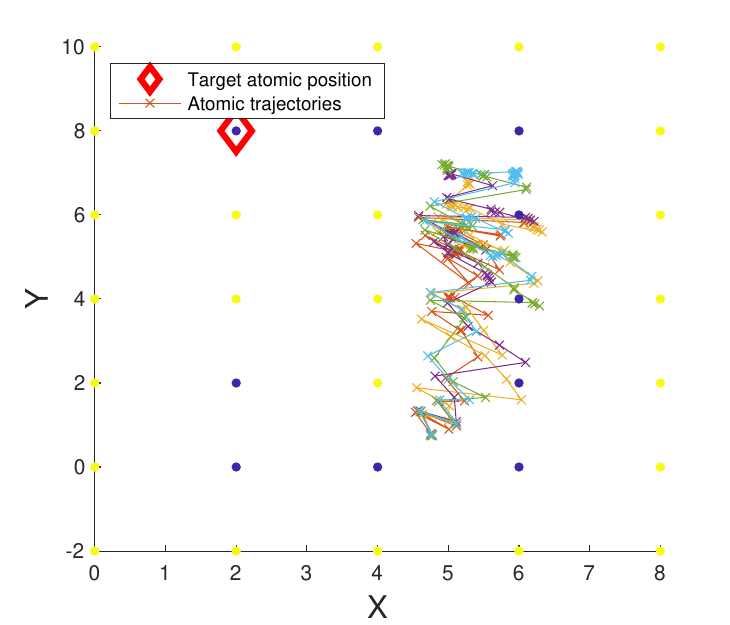}}\\
\subfloat[$10^5$ training episodes in Q-Learning]{\includegraphics[scale = 0.55]{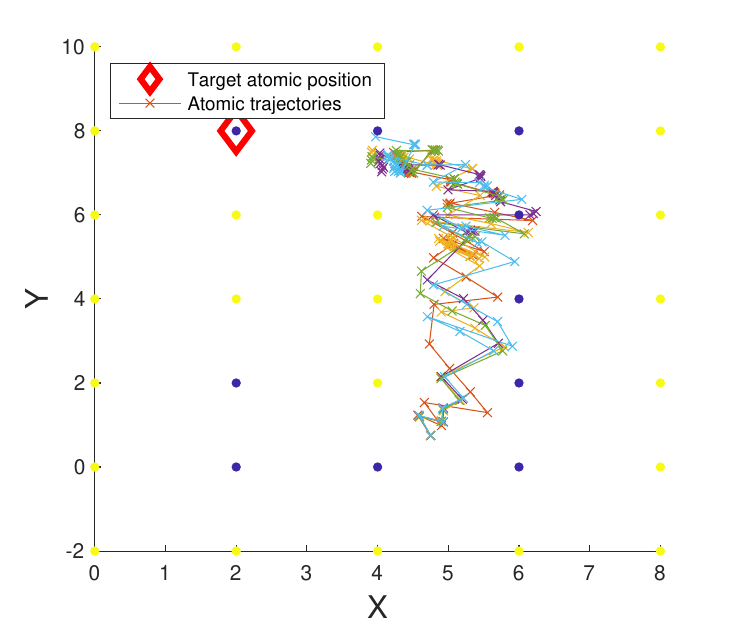}}
\subfloat[$10^6$ training episodes in Q-Learning]{\includegraphics[scale = 0.55]{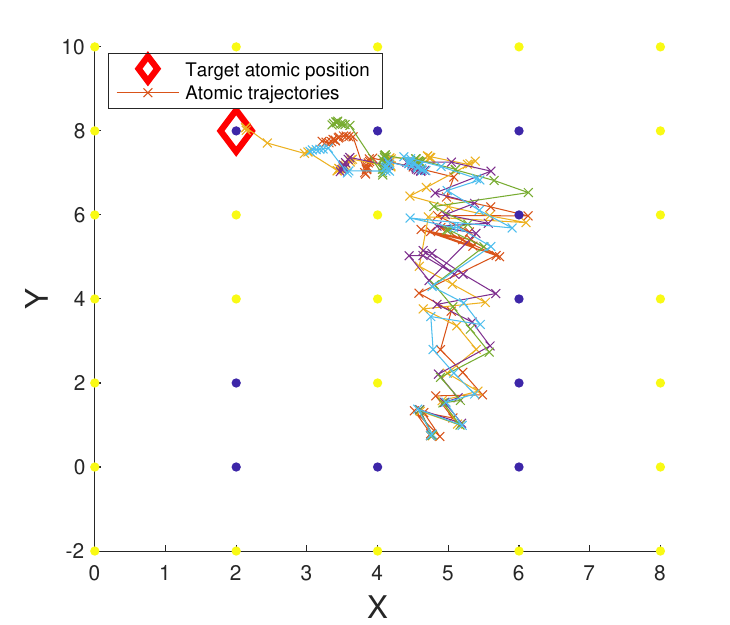}}
\end{center}
%\vspace{-0.75em}
\caption{Example 2. Performance of Q-Learning. }\label{Atomic_QL}
\end{figure}

In Figure \ref{Atomic_QL}, we show $5$ agent performance trajectories following the policy determined by the trained q-table with $10^3$, $10^4$, $10^5$, and $10^6$ training episodes in subplots (a), (b), (c), and (d), respectively. 
We can see from this Q-Learning performance experiment that with $10^3$ training episodes in the Q-Learning method, the agent could avoid the shallow potential ``trap atoms'' on the left side, which looks more promising at beginning but has some deep potential-well atoms as the ``barrier'' towards the destination atom (marked by the red diamond). However, with only $10^3$ episodes, the Q-Learning method cannot generate a sufficiently trained q-table that guides the agent towards the final destination. We can also observe from Figure \ref{Atomic_QL} that with longer and longer training procedures, the agent performance becomes better and better. With $10^6$ training episodes, one of those $5$ agents can finally arrive at $X_{destination}$, and the other $4$ agents also get close to the destination.

To show the advantageous performance of our algorithm, we also solve the RL problem for the atomic-forge technique by using our BAL method. To compare with the Q-Learning method, we use the same temporal partition, and we introduce the same uniform spatial partition with step-size $\Delta x = 0.1$ to approximate the state region $\mathcal{D}_{t_n}$ at each time $t_n$. Therefore, we approximate the continuous RL problem  Eq. \eqref{Atom-state}-\eqref{Atom-cost} with the same level of discretization accuracy as the Q-Learning method when using the BAL method.  To explore the environment, we use $Q = 100$ particles in the direct filter to estimate the environment parameter, and we carry out $K=10^3$ episodes in the training procedure.  The initial guess for the potential depth parameter of each atom is chosen as $\lambda_{guess} = 1$.
\begin{figure}[h!]
\begin{center}
\includegraphics[scale = 0.7]{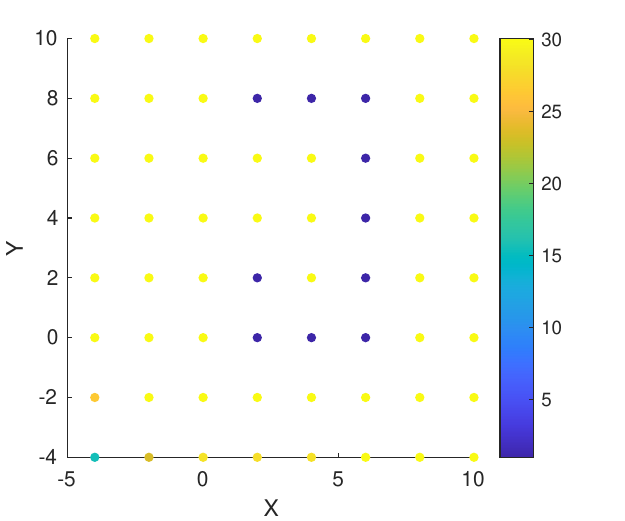}
%\vspace{-0.25em}
\end{center}
%\vspace{-0.75em}
\caption{Example 2. The estimated potential parameters on the material surface. }\label{Atomic_Estimated_Parameter} %\vspace{-0.5em}
\end{figure}
We first present the performance of direct filter based exploration in Figure \ref{Atomic_Estimated_Parameter}, where we also use the same color bar as we used in Figure \ref{Atomic_True_Parameter} to demonstrate the estimated potential depth parameter values corresponding to different background atoms. By comparing the estimated environment in Figure \ref{Atomic_Estimated_Parameter} with the true atomic environment in Figure \ref{Atomic_True_Parameter}, we can see that the environment learned in the BAL method is very similar to the true environment except for the left-bottom corner.  Note that the agent (i.e. the target atom) tries to move around the atoms with shallow potential-well (blue atoms) to save energy, and the agent does not have much experience near the left-bottom corner, which makes the left-bottom corner insufficiently explored.

\begin{figure}[h!]
\begin{center}
\includegraphics[scale = 0.7]{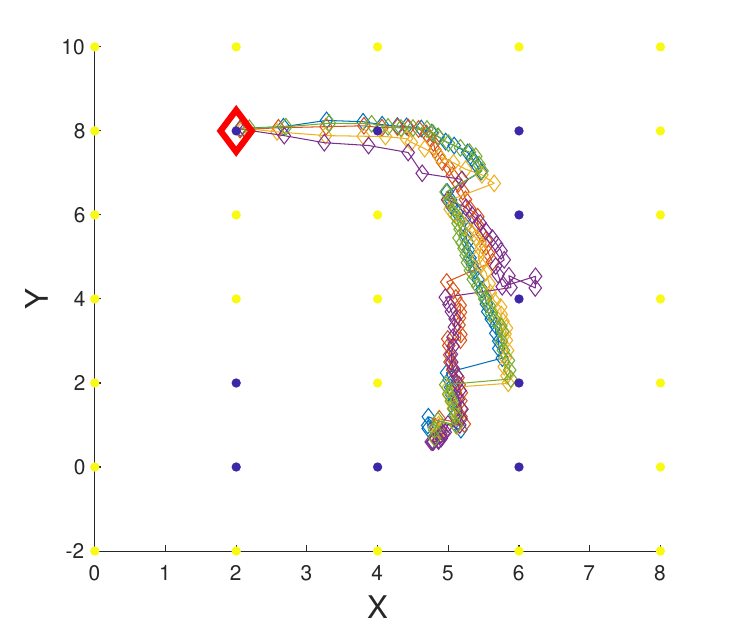}
%\vspace{-0.25em}
\end{center}
%\vspace{-0.75em}
\caption{Example 2. Performance of BAL.  }\label{Atomic_BAL} %\vspace{-0.5em}
\end{figure}
In Figure \ref{Atomic_BAL}, we present $5$ performance trajectories of the agent guided by the trained policy, which is obtained by the BAL method with $10^3$ training episodes. From this figure, we can see that all $5$ agents arrived at the designated target location at the terminal time. We also want to mention that due to the exploration mechanism in the BAL approach for RL, the agent has good understanding of the potential-well depths for the atoms near its route towards the destination. On the other hand, since the agent does not need to move around near the atoms at the left-bottom in the graph, it cannot learn the potential-well depth parameters very well for those atoms.

\begin{figure}[h!] %\vspace{-0.75em}
\begin{center}
\includegraphics[scale = 0.5]{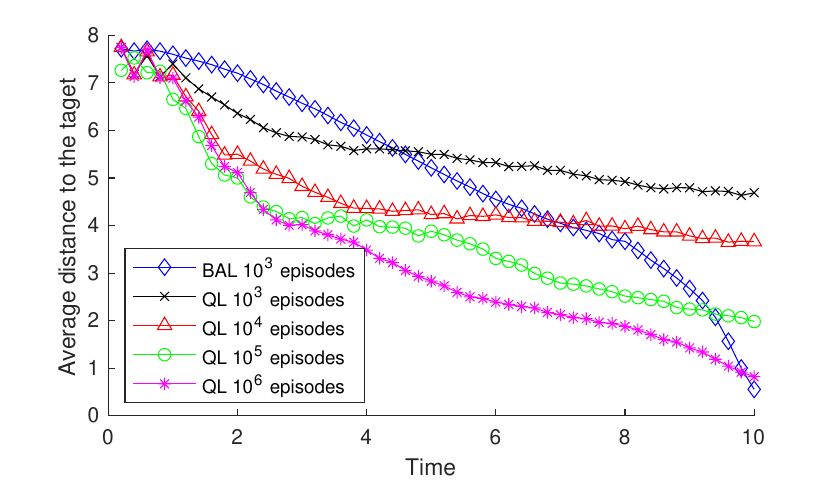}  
\end{center}
%\vspace{-0.75em}
\caption{Example 2. Comparison of average distance to the destination. }\label{Atomic_RMSE_Distance}
\end{figure}
To further demonstrate the performance comparison between the Q-Learning method and the BAL method, we calculate the average distance to the destination by letting the agent repeatedly run $100$ trajectories, and we present the average distance in Figure \ref{Atomic_RMSE_Distance}.  We can see from this figure that the BAL method with $10^3$ training episodes outperforms the Q-Learning method with $10^3$, $10^4$, and $10^5$ training episodes in terms of the average distance to the destination, and it slightly outperforms the Q-Learning method with $10^6$ training episodes with a very small margin. 
\begin{figure}[h!] %\vspace{-0.75em}
\begin{center}
\includegraphics[scale = 0.5]{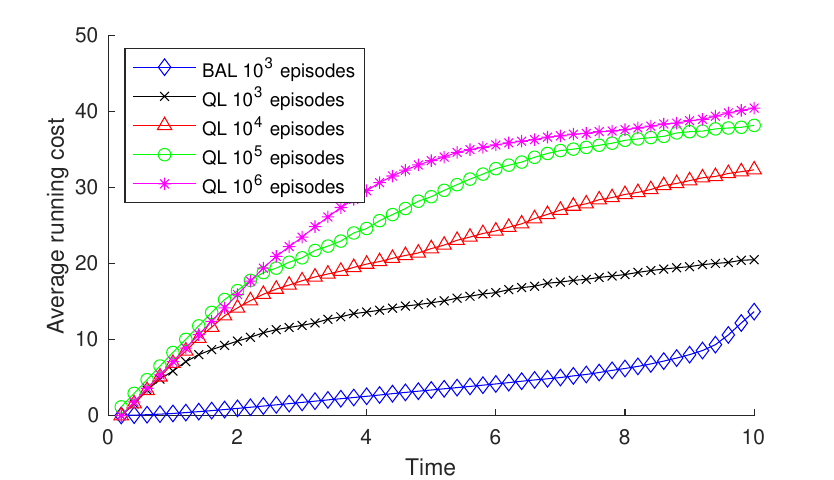}
%\subfloat[Average running costs]{\includegraphics[scale = 0.5]{./figs/Atomic_RMSE_Running}}
%\subfloat[Average overall costs]{\includegraphics[scale = 0.5]{./figs/Atomic_RMSE_Overall} } 
\end{center}
%\vspace{-0.75em}
\caption{Example 2. Comparison of average running costs. }\label{Atomic_RMSE}
\end{figure}
However, the criteria for the general performance of the agent is not only the distance to the designated destination. The consumption of energy for moving the agent, i.e. the running cost $\int_0^T |f_t| dt$, should also count for the performance. In Figure \ref{Atomic_RMSE}, we present the average running cost of each implementation with respect to time. We can see from this figure that the BAL requires much lower energy consumption compared with all the Q-Learning implementations -- regardless of the distance to the destination.

\subsection{Example 3: Reinforcement learning for continuous maneuvering of a robot in a maze.}\label{Maze}
In this example, we solve an RL problem for continuous maneuvering of a robot in a 2-dimensional maze, which is a continuous version of the maze solver problem \cite{2021maze, LSTM_Maze_01}.  We want to use this example to show the major advantage of the SMP approach over the DPP approach. Specifically, as an SMP type approach the BAL method treats the entire control procedure as a whole task, which allows the RL algorithm to design the control policy based on the predicted future trajectories of the agent. On the other hand, the Q-Learning method typically balances between the short-term optimal performance and the possible ultimate goal of the control task, and the design of the policy in the Q-Learning method does not rely on the comprehensive understanding of the entire environment. As a result, the Q-Learning designed policy could make the agent stuck in a local dilemma.

To proceed, we consider the following stochastic dynamics that describe the agent state
\begin{equation}\label{Robot-State}
\begin{aligned}
dX^{(1)}_t =& v \cos(\theta) dt + \sigma dW^{(1)}_t, \\
dX^{(2)}_t =& v \sin(\theta) dt + \sigma dW^{(2)}_t,
\end{aligned}
\end{equation}
where $\bm{X}_t = (X^{(1)}_t, X^{(2)}_t)^{\top}$ is the 2-dimensional location of a robot, the policy term $\bm{a} = (v, \theta)$ controls the velocity $v$ and the steering action $\theta$. The goal of the RL task in this example is to let the robot, i.e. the agent, learn how to arrive at a pre-designated destination $\bm{X}^{\ast}_T$ at the given terminal time $T$ with the lowest cost, and the cost function that we aim to minimize during the learning procedure is 
\begin{equation}\label{Robot-Cost}
J(\bm{a}) = \E[\int_0^T \lambda_x \|\bm{X}_t - \bm{X}_0\|_2^2 dt + F \|\bm{X}_T - \bm{X}^{\ast}_T\|_2^2],
\end{equation}
where $F = 20$ is a terminal cost constant that defines the amount of penalty for not being able to arrive at the destination $\bm{X}^{\ast}_T$ at the terminal time $T$, and $\int_0^T \lambda_x \|\bm{X}_t - \bm{X}_0\|_2^2 dt$ is the running cost term with an unknown parameter $\lambda_x$. Different from the RL problem for the atomic forge technique, in which we try to learn environment parameters that determine the state dynamics. In this work, as a RL problem for maze solver, the parameters that represent the unknowns in the environment are in the cost function $J(\bm{a})$ defined in Eq. \eqref{Robot-Cost}. Moreover, the RL problem that we try to solve in this example is different from a standard maze solver RL problem, in which the environment contains discrete obstacles (that add constant penalties to the cost if the agent touches them) and barriers (that stop the robot from getting through). In the RL problem Eq. \eqref{Robot-State} - \eqref{Robot-Cost}, we introduce a space-dependent cost parameter $\lambda_x$, which does not stop the agent from getting through, and the cost parameter $\lambda_x$ will \textit{continuously} add different levels of cost as the agent moving in the environment. Therefore, stepping into a small region with high cost parameter may not bring very high cost to the overall cost function $J(\bm{a})$. On the other hand, if the agent moves in a region with a fixed cost parameter, the faster the agent moves (or the farther the agent travels in a unit time) the more cost will be generated due to the accumulated running cost with respect to the distance that the agent traveled, i.e. $\|\bm{X}_t - \bm{X}_0\|_2$ in Eq. \eqref{Robot-Cost}.

To carry out numerical experiments, we introduce a temporal partition over time interval $[0, T]$ with $T = 20$, and we choose the time step-size $\Delta t = 0.2$, i.e. $N_T = 100$. The diffusion coefficient in Eq. \eqref{Robot-State} is chosen as $\sigma = 0.05$, which can bring relatively large amount of uncertainty to the state process given the length of the time interval and the time step-size. Also, we introduce a spatial partition to the environment by letting the x-dimension step-size $\Delta x = 0.2$ and y-dimension step-size $\Delta y = 0.25$, and we introduce a small base running cost $\lambda^{\text{base}} = 0.02$ everywhere in the state-space. The initial state of the agent is chosen as $\bm{X}_0 = (5, 4)^{\top}$, and the destination location is $\bm{X}^{\ast}_T = (5, 25)^{\top}$. 
\vspace{0.5em}

We first solve the above maze problem by using the Q-Learning method. Apparently, if the parameter $\lambda_x $ remains as the small invariant base cost $\lambda^{\text{base}}$ over the entire state space, the Q-Learning method would quickly converge and provide a policy that guides the agent to arrive at the destination directly. However, when an unknown high cost obstacle region appears, the maze problem could be more challenging. 
\begin{figure}[h!]
\begin{center}
\includegraphics[scale = 0.65]{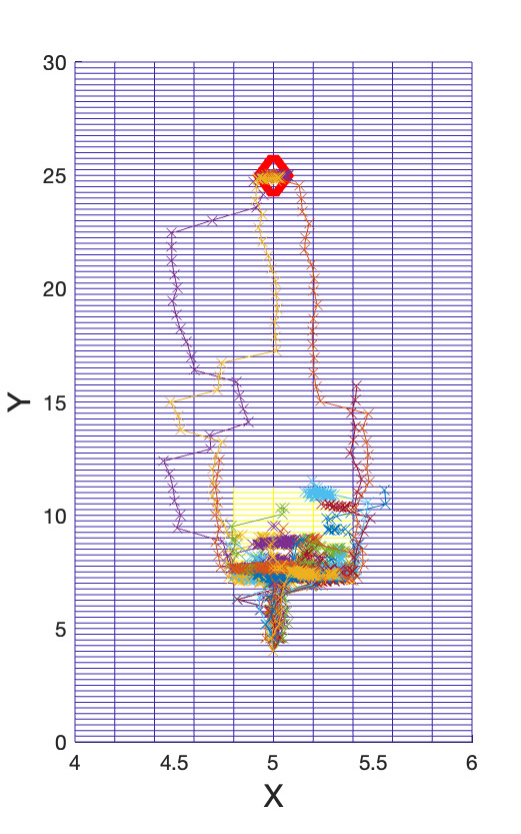}
%\vspace{-0.25em}
\end{center}
%\vspace{-0.75em}
\caption{Example 3. Performance of Q-Learning with a smaller obstacle region in the environment. }\label{Robot_QL_Small} %\vspace{-0.5em}
\end{figure}
In the first Q-Learning experiment, we put a rectangle region $[4.8, 5.4] \times [7, 11]$  in the environment and let the running cost coefficient in the region be $\lambda_x = 20$. Then, we train the q-table for $5 \times 10^5$ episodes and present $30$  testing agent trajectories using the trained q-table in Figure \ref{Robot_QL_Small}, where the red diamond shows the location of the destination $\bm{X}_T^{\ast}$, and the background mesh reflects the spatial partition for the state space. We can see from this figure that all the agents know that they should move up towards the destination. Once they step into the obstacle region, since they are still getting closer to the destination, the cost (or penalty) in the near future may not be large enough to stop them from moving forward. In other words, the high penalty of not being able to reach the destination may persuade the agent to overcome the short-term difficulties. However, when the agent accumulates large enough cost as it gets deeper into the obstacle region, the cost would only increase no matter where the agent goes. Therefore, we can see from Figure \ref{Robot_QL_Small} that most agents stop somewhere in the obstacle region. At the same time, we can also see from the figure that there are still several policy paths that would guide the agent to avoid the obstacle region, and it's more likely that those agents can arrive at the designated destination.

\begin{figure}[h!] %\vspace{-0.75em}
\begin{center}
\subfloat[Performance of Q-Learning: 30 agent trajectories]{\includegraphics[scale = 0.65]{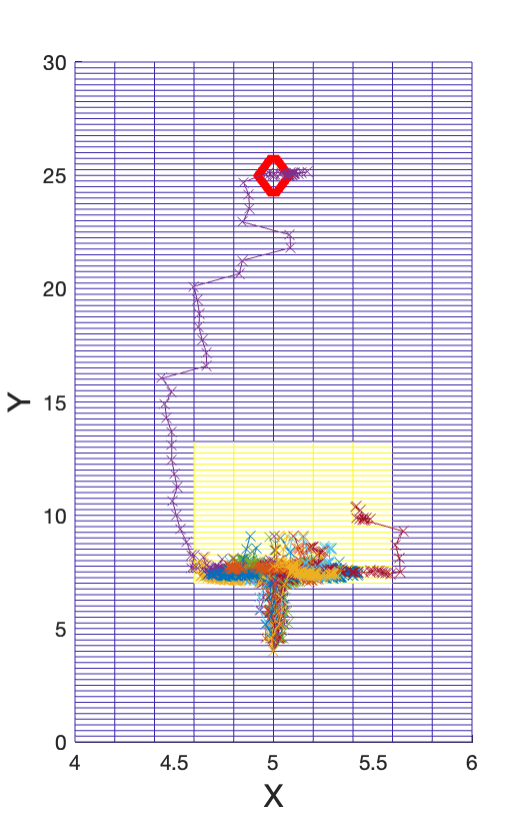} } \qquad
\subfloat[Performance of Q-Learning: 100 agent trajectories]{\includegraphics[scale = 0.65]{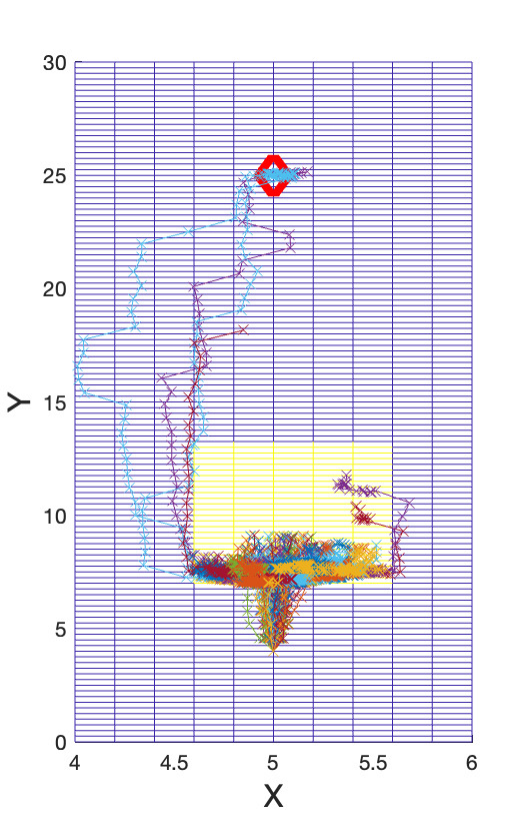} } 
\end{center}
%\vspace{-0.75em}
\caption{Example 2. Performance of Q-Learning with a larger obstacle region in the environment. }\label{Robot_QL_Large}
\end{figure}
To explore more along this direction and study a more challenging scenario for the Q-Learning method, we increase the size of the obstacle region to $[4.6, 5.6] \times [7, 13]$ and still train the q-table for $5 \times 10^5$ episodes. The performance of $30$ testing agent trajectories is plotted in Figure \ref{Robot_QL_Large} (a). From this figure, we can see more clearly how the agents moved horizontally and tried to move out of the obstacle region. One agent trajectory on the right hand side actually moved out of the obstacle region. Unfortunately, as this agent turned back towards the destination, it's trapped by the obstacle region again. The only agent that successfully arrived at the destination is plotted by the trajectory on the left. To further demonstrate the performance of the agent in this experiment, we present $100$ agent testing trajectories in Figure \ref{Robot_QL_Large} (b). From this figure, we can see that the ``successful'' policy followed by the trained q-table is on the left, and only the agents that follow the left-side-policy could arrive at the destination. 

%From this experiment, we can also conjecture that by training the agent with more episodes, the Q-table may eventually trained enough to create more ``successful'' policies. However, in reality, the environment could be much more than what we see in Figure \ref{Robot_QL_Small} and Figure \ref{Robot_QL_Large}.

\begin{figure}[h!]
\begin{center}
\includegraphics[scale = 0.85]{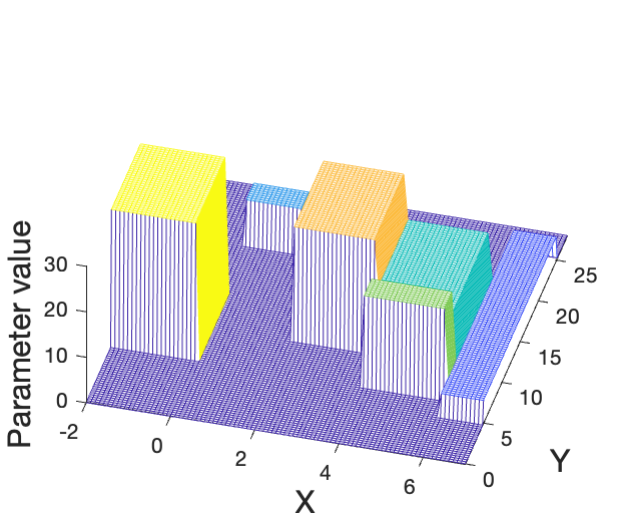}
%\vspace{-0.25em}
\end{center}
%\vspace{-0.75em}
\caption{Example 3. The 3D view of the true environment of the maze. The heights of the obstacle regions show different cost parameter values. }\label{Robot_True_Environment} %\vspace{-0.5em}
\end{figure}

In Figure \ref{Robot_True_Environment}, we design a much more complicated environment with $6$ obstacle regions, and each region has a different cost parameter value. Specifically, from the right to the left we let $\lambda_x^1 = 5$, $\lambda_x^2 = 20$, $\lambda_x^3 = 15$, $\lambda_x^4 = 25$, $\lambda_x^5 = 10$, $\lambda_x^6 = 30$. 
The different running cost parameter values are presented by the heights of the obstacle regions in Figure \ref{Robot_True_Environment}.
\begin{figure}[h!]
\begin{center}
\includegraphics[scale = 0.75]{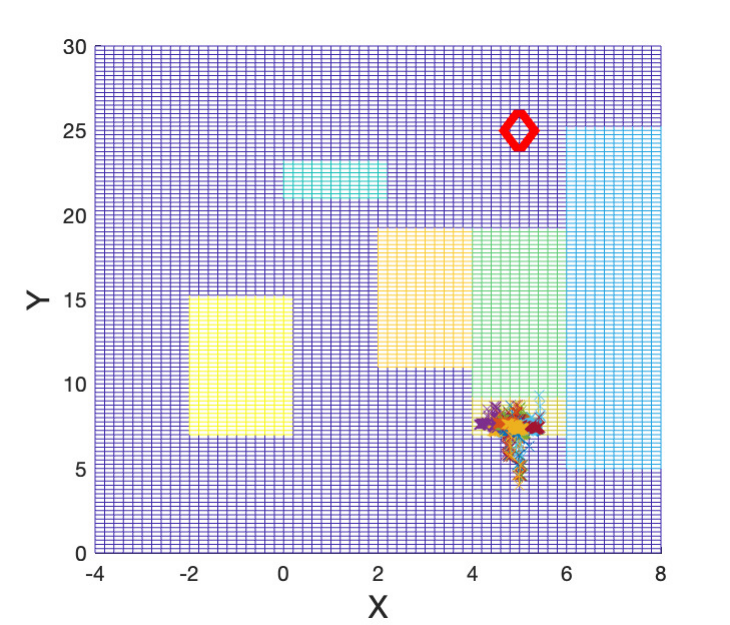}
%\vspace{-0.25em}
\end{center}
%\vspace{-0.75em}
\caption{Example 3. Performance of Q-Learning in the maze presented in Figure \ref{Robot_True_Environment}.  }\label{Robot_QL_Full} %\vspace{-0.5em}
\end{figure}
In Figure \ref{Robot_QL_Full}, we show the performance of the agent with $30$ testing trajectories that follow the policy obtained by the Q-Learning method with $5 \times 10^5$ training episodes. There's no surprise that the agent cannot find a path to avoid all the obstacles and arrive at the destination. Although we can conjecture that by training the agent with more and more episodes, the q-table may eventually be well-trained enough to create some ``successful'' policies. However, for such a complicated environment with so many trapping obstacle regions, it would be very computationally expensive for Q-Learning to find a path to arrive at the destination.

\vspace{0.5em}

In the following experiments, we show the success of the BAL method in solving the RL problem Eq. \eqref{Robot-State} - \eqref{Robot-Cost}. To implement the BAL method, we use $Q = 100$ particles to carry out the direct filter based parameter estimation for exploration, and we carry out $K=1000$ training episodes in the BAL algorithm. The temporal partition and the spatial partition that we use for the BAL method are the same as the Q-Learning method. In this experiment, we don't assume that the agent knows there are altogether $6$ obstacle regions. Instead, we use the direct filter based exploration technique to estimate the running cost parameter in every artificially partitioned spatial block with partition size $\Delta x \times \Delta y$. 
\begin{figure}[h!]
\begin{center}
\includegraphics[scale = 0.85]{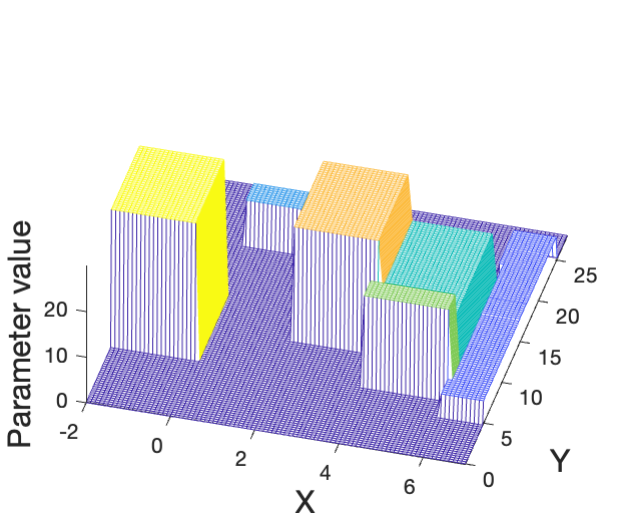}
%\vspace{-0.25em}
\end{center}
%\vspace{-0.75em}
\caption{Example 3. The estimated estimated environment of the maze obtained by the BAL method. }\label{Robot_Estimated_Environment} %\vspace{-0.5em}
\end{figure}

In Figure \ref{Robot_Estimated_Environment}, we first present the estimated environment learned by the BAL method, where the background mesh shows all the artificially partitioned spatial  blocks that determine the size of environment parameter. We can see by comparing Figure  \ref{Robot_Estimated_Environment} with Figure \ref{Robot_True_Environment} that the BAL method successfully recovered the environment, and the estimate for every obstacle region is very accurate -- in both the size of each obstacle region and the value of each cost parameter. Although we only implemented $1000$ training episodes, since the direct filter updates the estimate for the environment at every time instant in each training episode, the parameter estimation algorithm has sufficiently incorporated the agent trial states into the exploration procedure, and this helps the agent have very good understanding of the environment.

\begin{figure}[h!]
\begin{center}
\includegraphics[scale = 0.75]{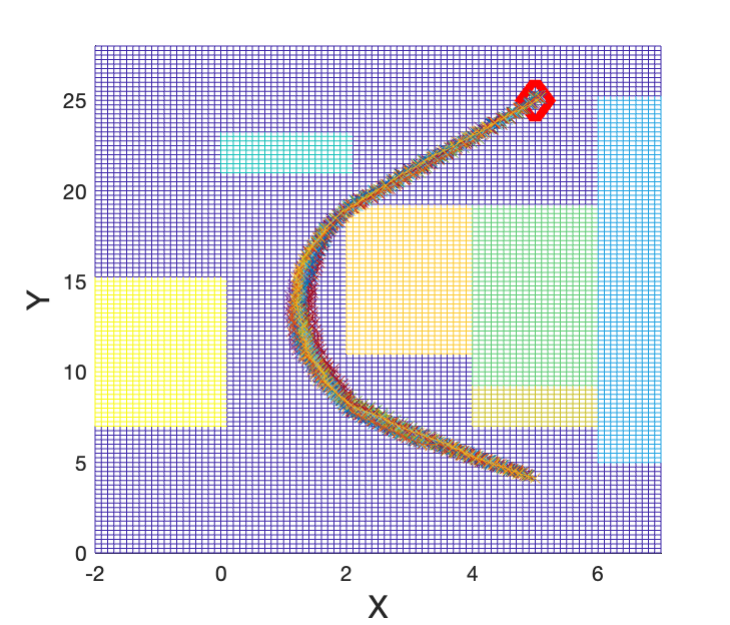}
%\vspace{-0.25em}
\end{center}
%\vspace{-0.75em}
\caption{Example 3. Performance of the BAL method in the maze.  }\label{Robot_BAL_Point} %\vspace{-0.5em}
\end{figure}
In Figure \ref{Robot_BAL_Point}, we present $30$ agent trajectories guided by the BAL trained policy with initial state $X_0 = (5, 4)^{\top}$ and terminal destination is still chosen as $\bm{X}^{\ast}_T = (5, 25)^{\top}$. We can see from this figure that all the agent trajectories follow very smooth paths towards the destination, and they all arrived at the destination at the terminal time.

\begin{figure}[h!]
\begin{center}
\includegraphics[scale = 0.75]{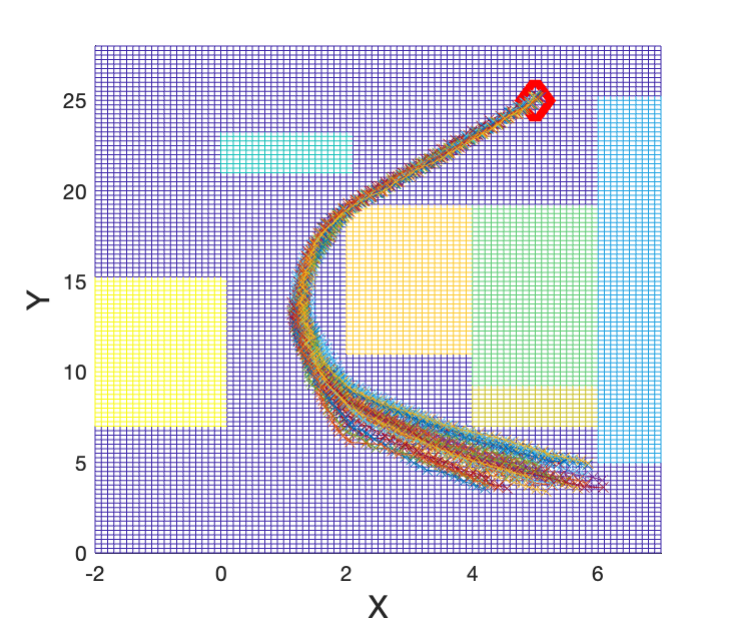}
%\vspace{-0.25em}
\end{center}
%\vspace{-0.75em}
\caption{Example 3. Performance of the BAL method in the maze with random initial states. }\label{Robot_BAL_Region} %\vspace{-0.5em}
\end{figure}

To further demonstrate the performance of the BAL method in this example and show the robustness of our method, we let the agent start from randomly picked initial states in the spatial area $[4, 6] \times [3, 5]$, and we present $30$ agent trajectories following the BAL trained policy in Figure \ref{Robot_BAL_Region}. 
We can see from this figure that no matter where the agent started, it can always find the right path towards the destination, and it can perfectly avoid the obstacle regions.

\vspace{0.5em}
By comparing the performance of the Q-Learning trained policy with the BAL trained policy from the above experiments, we can see that the BAL method clearly outperforms the Q-Learning method in this continuous maze solver RL problem.

\section{Conclusions}\label{Conclusions}
In this work, we developed a stochastic maximum principle (SMP) approach for solving the reinforcement learning (RL) problem in the case that the environment can be parameterized. To explore the environment, we introduced a direct filter method as an online parameter estimation method to learn the environment parameters during the training procedure, and the exploitation task is carried out through an efficient backward action learning (BAL) algorithm for finding the optimal policy under the SMP framework. The main advantage of such an SMP approach, compared with dynamic programming principle (DPP) based methods, is that the gradient of the cost with respect to the control aims to find the improvement direction for the control over the entire performance period, which could potentially provide better overall performance for the long-run. In contrast, the standard temporal difference (TD) learning methods under the DPP framework, such like Q-Learning methods, only consider short-term rewards or penalties, which could ignore possible future opportunities. Especially, like the numerical experiments presented in Example \ref{Maze}, if the environment is carefully designed so that the agent is deeply trapped in a situation where no short-term strategy can lead it out, any TD learning method will struggle unless it can consider long-term predictions. However, TD learning with long-term predictions can be very difficult when the environment is unknown. Although our BAL method also needs to consider long-term predictions, the direct filter method can be effectively combined with the BAL algorithm to establish a comprehensive understanding of the environment, and this allows us to dynamically learn the entire environment while searching for the long-term optimal policy.  The major drawback of our BAL method for solving the RL problem is that we require a parameterization for the environment. But for applications of RL in science, such like the atomic forge technique, physics knowledge and pre-defined models are necessary to make the RL problem ``physics-informed''. In this case, having a parameterized environment is a very reasonable assumption.

\bibliographystyle{plain}

\begin{thebibliography}{10}

\bibitem{MCMC-PF}
C.~Andrieu, A.~Doucet, and R.~Holenstein.
\newblock Particle markov chain monte carlo methods.
\newblock {\em J. R. Statist. Soc. B}, 72(3):269--342, 2010.

\bibitem{Bao_parameter}
R.~Archibald, F.~Bao, and X.~Tu.
\newblock A direct filter method for parameter estimation.
\newblock {\em J. Comput. Phys.}, 398:108871, 17, 2019.

\bibitem{Bao_Control_20}
R.~Archibald, F.~Bao, J.~Yong, and T.~Zhou.
\newblock An efficient numerical algorithm for solving data driven feedback
  control problems.
\newblock {\em Journal of Scientific Computing}, 85(51), 2020.

\bibitem{Bao_SNN_22}
Richard Archibald, Feng Bao, Yanzhao Cao, and Hui Sun.
\newblock Convergence analysis for training stochastic neural networks via
  stochastic gradient descent, 2022.

\bibitem{Bao_EAJAM20}
Richard Archibald, Feng Bao, and Jiongmin Yong.
\newblock A stochastic gradient descent approach for stochastic optimal
  control.
\newblock {\em East Asian Journal on Applied Mathematics}, 10(4):635--658,
  2020.

\bibitem{LSTM_Maze_01}
Bram Bakker.
\newblock Reinforcement learning with long short-term memory.
\newblock {\em Advances in neural information processing systems}, 14, 2001.

\bibitem{Bao_IJUQ19}
F.~Bao, Y.~Cao, and H.~Chi.
\newblock Adjoint forward backward stochastic differential equations driven by
  jump processes and its application to nonlinear filtering problems.
\newblock {\em International Journal of Uncertainty Quantification},
  9(2):143--159, 2019.

\bibitem{Bao_Cogan20}
F.~Bao, N.~Cogan, A.~Dobreva, and R.~Paus.
\newblock Data assimilation of synthetic data as a novel strategy for
  predicting disease progression in alopecia areata.
\newblock {\em Mathematical Medicine and Biology: A Journal of the IMA}, 2021.

\bibitem{Bao_first}
Feng Bao, Yanzhao Cao, Amnon Meir, and Weidong Zhao.
\newblock A first order scheme for backward doubly stochastic differential
  equations.
\newblock {\em SIAM/ASA J. Uncertain. Quantif.}, 4(1):413--445, 2016.

\bibitem{BSDE_filter}
Feng Bao and Vasileios Maroulas.
\newblock Adaptive meshfree backward {SDE} filter.
\newblock {\em SIAM J. Sci. Comput.}, 39(6):A2664--A2683, 2017.

\bibitem{cd2002}
D.~Crisan and A.~Doucet.
\newblock A survey of convergence results on particle filtering methods for
  practitioners.
\newblock {\em IEEE Trans. Sig. Proc.}, 50(3):736--746, 2002.

\bibitem{Bao_Atom20}
O.~Dyck, M.~Ziatdinov, S.~Jesse, F.~Bao, A.~Yousefzadi~Nobakht, A.~Maksov, B.G.
  Sumpter, R.~Archibald, K.J.H. Law, and S.V. Kalinin.
\newblock Probing potential energy landscapes via electron-beam-induced single
  atom dynamics.
\newblock {\em Acta Materialia}, 203:116508, 2021.

\bibitem{Bao_Carbon_22}
Ondrej Dyck, Feng Bao, Maxim Ziatdinov, Ali~Yousefzadi Nobakht, Kody Law, Artem
  Maksov, Bobby~G. Sumpter, Richard Archibald, Stephen Jesse, Sergei~V.
  Kalinin, and David~B. Lingerfelt.
\newblock Strain-induced asymmetry and on-site dynamics of silicon defects in
  graphene.
\newblock {\em Carbon Trends}, 9:100189, 2022.

\bibitem{Q_Learning_1997fuzzy}
Pierre~Yves Glorennec and Lionel Jouffe.
\newblock Fuzzy q-learning.
\newblock In {\em Proceedings of 6th international fuzzy systems conference},
  volume~2, pages 659--662. IEEE, 1997.

\bibitem{GPM_2017}
Bo~Gong, Wenbin Liu, Tao Tang, Weidong Zhao, and Tao Zhou.
\newblock An efficient gradient projection method for stochastic optimal
  control problems.
\newblock {\em SIAM J. Numer. Anal.}, 55(6):2982--3005, 2017.

\bibitem{particle-filter}
N.J Gordon, D.J Salmond, and A.F.M. Smith.
\newblock Novel approach to nonlinear/non-gaussian bayesian state estimation.
\newblock {\em IEE PROCEEDING-F}, 140(2):107--113, 1993.

\bibitem{Q_learning_continuous_16}
Shixiang Gu, Timothy Lillicrap, Ilya Sutskever, and Sergey Levine.
\newblock Continuous deep q-learning with model-based acceleration.
\newblock In {\em International conference on machine learning}, pages
  2829--2838. PMLR, 2016.

\bibitem{Kalinin-Atom}
S.~Kalinin, A.~Borisevich, and S.~Jesse.
\newblock Fire up the atom forge.
\newblock {\em Nature}, 22 November 2016.

\bibitem{Kang-PF}
Kai Kang, Vasileios Maroulas, Ioannis Schizas, and Feng Bao.
\newblock Improved distributed particle filters for tracking in a wireless
  sensor network.
\newblock {\em Comput. Statist. Data Anal.}, 117:90--108, 2018.

\bibitem{SDE1}
P.~E. Kloeden and E.~Platen.
\newblock {\em Numerical solution of stochastic differential equations},
  volume~23 of {\em Applications of Mathematics (New York)}.
\newblock Springer-Verlag, Berlin, 1992.

\bibitem{ED_RL_21}
Viraj Mehta, Biswajit Paria, Jeff Schneider, Stefano Ermon, and Willie
  Neiswanger.
\newblock An experimental design perspective on model-based reinforcement
  learning, 2021.

\bibitem{MTAC2012}
M.~Morzfeld, X.~Tu, E.~Atkins, and A.~J. Chorin.
\newblock A random map implementation of implicit filters.
\newblock {\em J. Comput. Phys.}, 231(4):2049--2066, 2012.

\bibitem{Ali-Atom}
Ali~Yousefzadi Nobakht, Ondrej Dyck, David~B Lingerfelt, Feng Bao, Maxim
  Ziatdinov, Artem Maksov, Bobby~G Sumpter, Richard Archibald, and Sergei
  V~Kalinin Stephen~Jesse, and Kody~JH Law.
\newblock Reconstruction of effective potential from statistical analysis of
  dynamic trajectories.
\newblock {\em AIP Advances}, 10:065034, 2020.

\bibitem{Q_learning_94}
Jing Peng and Ronald~J Williams.
\newblock Incremental multi-step q-learning.
\newblock In {\em Machine Learning Proceedings 1994}, pages 226--232. Elsevier,
  1994.

\bibitem{Peng_control}
Shi~Ge Peng.
\newblock A general stochastic maximum principle for optimal control problems.
\newblock {\em SIAM J. Control Optim.}, 28(4):966--979, 1990.

\bibitem{Sutton_RL_2nd}
Richard S.~Sutton and Andrew G.~Barto.
\newblock Reinforcement learning: An introduction second edition: 2014, 2015.
\newblock 2014.

\bibitem{TD_Learning_95}
Gerald Tesauro et~al.
\newblock Temporal difference learning and td-gammon.
\newblock {\em Communications of the ACM}, 38(3):58--68, 1995.

\bibitem{epsilon_greedy11}
Michel Tokic and G{\"u}nther Palm.
\newblock Value-difference based exploration: adaptive control between
  epsilon-greedy and softmax.
\newblock In {\em Annual conference on artificial intelligence}, pages
  335--346. Springer, 2011.

\bibitem{DRL_ED_2022}
Neythen~J. Treloar, Nathan Braniff, Brian Ingalls, and Chris~P. Barnes.
\newblock Deep reinforcement learning for optimal experimental design in
  biology.
\newblock {\em bioRxiv}, 2022.

\bibitem{DRL_16}
Hado Van~Hasselt, Arthur Guez, and David Silver.
\newblock Deep reinforcement learning with double q-learning.
\newblock In {\em Proceedings of the AAAI conference on artificial
  intelligence}, volume~30, 2016.

\bibitem{XYZ_RL_20}
Haoran Wang, Thaleia Zariphopoulou, and Xun~Yu Zhou.
\newblock Reinforcement learning in continuous time and space: A stochastic
  control approach.
\newblock {\em Journal of Machine Learning Research}, 21(198):1--34, 2020.

\bibitem{Q_learning_92}
Christopher~JCH Watkins and Peter Dayan.
\newblock Q-learning.
\newblock {\em Machine learning}, 8(3):279--292, 1992.

\bibitem{Yong_control}
Jiongmin Yong and Xun~Yu Zhou.
\newblock {\em Stochastic controls}, volume~43 of {\em Applications of
  Mathematics (New York)}.
\newblock Springer-Verlag, New York, 1999.
\newblock Hamiltonian systems and HJB equations.

\bibitem{ZhangJ_BSDE}
Jianfeng Zhang.
\newblock A numerical scheme for {BSDE}s.
\newblock {\em Ann. Appl. Probab.}, 14(1):459--488, 2004.

\bibitem{2021maze}
Xiaoping Zhang, Yihao Liu, Dunli Hu, and Lei Liu.
\newblock A maze robot autonomous navigation method based on curiosity and
  reinforcement learning.
\newblock In {\em The 7th Int. Workshop on Advanced Computational Intelligence
  and Intelligent Informatics (IWACIII 2021), Article}, number M1-6, page~1,
  2021.

\bibitem{Zhao_multi}
Weidong Zhao, Yu~Fu, and Tao Zhou.
\newblock New kinds of high-order multistep schemes for coupled forward
  backward stochastic differential equations.
\newblock {\em SIAM J. Sci. Comput.}, 36(4):A1731--A1751, 2014.

\end{thebibliography}

\end{document}